\documentclass[a4paper,11pt]{article}
\usepackage[latin1]{inputenc}
\usepackage[francais,english]{babel}    % Pour du fran\c{c}ais
\usepackage{amssymb}
\usepackage{amsmath}
\usepackage{amsthm}
\usepackage{multicol}
\usepackage{graphicx}
\usepackage{color}
\usepackage[T1]{fontenc}
\usepackage{yfonts}
\usepackage{placeins}
\usepackage{array}
\usepackage{dsfont}
\usepackage{stmaryrd}
\usepackage{verbatim}
\usepackage{caption}
\usepackage{transparent}
\usepackage{bbold}
\usepackage{hyperref}
\usepackage{enumitem}
\usepackage{tikz}
\usetikzlibrary{patterns}
\usepackage{tikz-3dplot}

\newtheorem{thm}{Theorem}[section]
\newcommand{\bthm}{\begin{thm}}
\newcommand{\ethm}{\end{thm}}

\newtheorem{thmi}{Theorem}
\newcommand{\bthmi}{\begin{thmi}}
\newcommand{\ethmi}{\end{thmi}}

\newtheorem{cori}[thmi]{Corollary}
\newcommand{\bcori}{\begin{cori}}
\newcommand{\ecori}{\end{cori}}

\newtheorem{mthm}{Theorem}
\newcommand{\bmthm}{\begin{mthm}}
\newcommand{\emthm}{\end{mthm}}

\newtheorem{mcor}[mthm]{Corollary}
\newcommand{\bmcor}{\begin{mcor}}
\newcommand{\emcor}{\end{mcor}}
\newtheorem{mconj}[mthm]{Conjecture}
\newcommand{\bmconj}{\begin{mconj}}
\newcommand{\emconj}{\end{mconj}}
\newtheorem{mpro}[mthm]{Proposition}
\newcommand{\bmpro}{\begin{mpro}}
\newcommand{\empro}{\end{mpro}}

\newtheorem*{conj}{Conjecture}
\newcommand{\bconj}{\begin{conj}}
\newcommand{\econj}{\end{conj}}

\newtheorem*{question}{Question}
\newcommand{\bq}{\begin{question}}
\newcommand{\eq}{\end{question}}

\newtheorem*{thn}{Theorem}
\newcommand{\bthn}{\begin{thn}}
\newcommand{\ethn}{\end{thn}}

\newtheorem{exo}{Exercise}
\newcommand{\bex}{\begin{exo}}
\newcommand{\eex}{\end{exo}}

\newtheorem{sol}{Solution}
\newcommand{\bsol}{\begin{sol}}
\newcommand{\esol}{\end{sol}}

\newtheorem{pro}[thm]{Proposition}
\newcommand{\bpro}{\begin{pro}}
\newcommand{\epro}{\end{pro}}

\newtheorem{cor}[thm]{Corollary}
\newcommand{\bcor}{\begin{cor}}
\newcommand{\ecor}{\end{cor}}

\newtheorem{lem}[thm]{Lemma}
\newcommand{\blem}{\begin{lem}}
\newcommand{\elem}{\end{lem}}

\theoremstyle{definition}

\newtheorem{defi}[thm]{Definition}
\newcommand{\bdf}{\begin{defi}}
\newcommand{\edf}{\end{defi}}

\newtheorem*{defis}{Definition}
\newcommand{\bdfs}{\begin{defis}}
\newcommand{\edfs}{\end{defis}}

\newtheorem*{rmk}{Remark}
\newcommand{\brk}{\begin{rmk} \upshape}
\newcommand{\erk}{\end{rmk}}

\newtheorem*{rmks}{Remarks}
\newcommand{\brks}{\begin{rmks} \upshape}
\newcommand{\erks}{\end{rmks}}

\newtheorem*{exe}{Example}
\newcommand{\bexe}{\begin{exe} \upshape}
\newcommand{\eexe}{\end{exe}}

\newtheorem*{exes}{Examples}
\newcommand{\bexes}{\begin{exes} \upshape}
\newcommand{\eexes}{\end{exes}}

\newtheorem*{pre}{Proof}
\newcommand{\bp}{\begin{pre} \upshape}
\newcommand{\ep}{\hfill \qed \end{pre}}
\newcommand{\epp}{\end{pre}}

\newcommand{\beq}{\begin{eqnarray*}}
\newcommand{\eeq}{\end{eqnarray*}}

\newcommand{\beqn}{\begin{equation}}
\newcommand{\eeqn}{\end{equation}}

\newcommand{\ben}{\begin{enumerate}}
\newcommand{\een}{\end{enumerate}}
\newcommand{\bit}{\begin{itemize} \renewcommand{\labelitemi}{$\bullet$} \renewcommand{\labelitemii}{$\star$}}
\newcommand{\eit}{\end{itemize}}

\newcommand{\bfg}{
\begin{figure}[H]
\begin{center}}
\newcommand{\efg}{
\end{center}
\end{figure}
\FloatBarrier}

\newcolumntype{M}[1]{>{\raggedright}m{#1}}

\newcommand{\R}{\mathbb{R}}

\newcommand{\N}{\mathbb{N}}
\newcommand{\Z}{\mathbb{Z}}

\newcommand{\K}{\mathbb{K}}
\renewcommand{\k}{\mathbb{k}}

\renewcommand{\tilde}{\widetilde}

\newcommand{\Ker}{\operatorname{Ker}}

\newcommand{\Mod}{\operatorname{Mod}}

\newcommand{\eps}{\varepsilon}
\newcommand{\st}{\, | \,}

\newcommand{\ra}{\rightarrow}

\renewcommand{\geq}{\geqslant}
\renewcommand{\leq}{\leqslant}

\newcommand{\GL}{\operatorname{GL}}
\newcommand{\SL}{\operatorname{SL}}

\newcommand{\PGL}{\operatorname{PGL}}

\newcommand{\<}{\langle}
\renewcommand{\>}{\rangle}

\newcommand{\mk}{\medskip}

\makeatletter
\def\Ddots{\mathinner{\mkern1mu\raise\p@
\vbox{\kern7\p@\hbox{.}}\mkern2mu
\raise4\p@\hbox{.}\mkern2mu\raise7\p@\hbox{.}\mkern1mu}}
\makeatother

\makeatletter
\def\maketitles{%
  \null
  \thispagestyle{empty}%
  \vfill
  \begin{center}\leavevmode
    \normalfont
    {\LARGE \@title\par}%
    \vskip 1.2cm
    {\large \@author\par}%
    \vskip 1.2cm
    {\large \@subtitle\par}%
    \vskip 0.8cm
    {\large \@date\par}%
  \end{center}%
  \vfill
  \null
  \cleardoublepage
  }

\def\date#1{\def\@date{#1}}
\def\author#1{\def\@author{#1}}
\def\title#1{\def\@title{#1}}
\def\subtitle#1{\def\@subtitle{#1}}
\makeatother

\usepackage{float}
\newcommand{\Cay}{\operatorname{Cay}}

\usetikzlibrary{arrows}

\begin{document}

\title{Lattices, Garside structures and weakly modular graphs}
\author{Thomas Haettel\thanks{Thomas Haettel, IMAG, Univ Montpellier, CNRS, France, and IRL 3457, CRM-CNRS, Universit\'{e} de Montr\'{e}al, Canada.}\, and Jingyin Huang\thanks{Jingyin Huang, Department of Mathematics, The Ohio State University, 100 Math Tower, 231 W 18th Ave, Columbus, OH 43210, U.S.}}
\date{\today}

\selectlanguage{english}

\maketitle

\begin{center}
\begin{minipage}{0.8\textwidth}
\textsc{Abstract.} In this article we study combinatorial non-positive curvature aspects of various simplicial complexes with natural $\widetilde A_n$ shaped simplicies, including Euclidean buildings of type $\widetilde A_n$ and Cayley graphs of Garside groups and their quotients by the Garside elements. All these examples fit into the more general setting of lattices with order-increasing $\Z$-actions and the associated lattice quotients proposed in a previous work by the first named author. We show that both the lattice quotients and the lattices themselves give rise to weakly modular graphs, which is a form of combinatorial non-positive curvature. We also show that several other objects fit into this setting of lattices/lattice quotients, including Artin complexes of Artin-Tits groups of type $\widetilde A_n$, a class of arc complexes and weak Garside groups arising from a categorical Garside structure in the sense of Bessis. Hence our result also implies to these objects and shows that they give weakly modular graphs.
Along the way, we also clarify the relationship between categorical Garside structures, lattices with $\Z$ action and different classes of complexes studied this article. We use this point of view to describe the first examples of Garside groups with exotic properties, like non-linearity or rigidity results.
\end{minipage}
\end{center}

\let\thefootnote\relax\footnotetext{{\bf Keywords} : Buildings, Artin-Tits groups, Garside groups, nonpositive curvature, weakly modular graphs. {\bf AMS codes} : 20E42, 20F36, 20F55, 05B35, 06A12, 20F65, 05C25}

\mk

The authors would like to acknowledge the very deep influence of Jacques Tits in many topics relating algebra and geometry, notably buildings and Artin-Tits groups, which are both at the core of this article.

\mk

The topic of combinatorial non-positive curvature (CNPC) lies in the intersection of metric graph theory and geometric group theory. We refer to \cite{chalopin_chepoi_hirai_osajda} for a detailed discussion of the context and motivation for CNPC. The basic idea is to identify local combinatorial patterns of graphs or complexes that lead to standard consequences of non-positive curvature, e.g. propagation of these combinatorial patterns from local to global in the spirit of the classical Cartan-Hadamard theorem, control of isoperimetric inequalities, existence of nice combings, fixed point properties, asphericity etc. Pioneering examples of CNPC include small cancellation theory, Gromov's flagness condition \cite{gromov_hyperbolic_groups}, and systolic complexes \cite{januszkiewicz2006simplicial}.
For a group, being able to act on graphs or complexes that satisfy some form of CNPC usually has strong implications on the structure of this group. Thus it is of great interest to construct such actions - there are many works in this direction, and we simply mention a few which are closely related to this article \cite{bestvina_artin,huang2020large,huang_osajda_helly,chalopin2020helly,munro2019weak,hoda:crystallographic,cumplido2020parabolic,haettel_injective_buildings,soergel2021systolic,haettel2021lattices,blufstein2021parabolic,haettel2021locally}.

%Among the various form of CNPC, we want to highlight the notion of \emph{weakly modular graph}. The definition requires balls in the graphs satisfies a weak form of convexity (see Definition???). Rather interestingly, these form of convexity propagates from local to global. Weakly modular graphs serve a common generalization to several other important forms of CNPC, including median, systolic, bucolic, Helly. Moreover, fundamental properties of these other forms of CNPC factor though fundamental properties of weakly modular graphs.

There is a strong connection between affine Coxeter groups and forms of CNPC. Each such Coxeter group gives a Euclidean polyhedron which is the fundamental domain of the action of this group on the associated Euclidean space. Then one can try to build more complicated spaces using such Euclidean polyhedra, and ask whether there is a specific local combinatorial pattern of assembling these Euclidean polyhedrons such that the resulting space is non-positively curved in some sense. For example, systolic complexes and bridged graphs \cite{chepoi2000graphs,januszkiewicz2006simplicial} describe a form of CNPC for spaces made of equilateral triangles; CAT(0) cube complexes and median graphs describe another form of CNPC for spaces made of unit cubes, quasi-medians graphs and bucolic complexes can be viewed as forms of CNPC for spaces made with prisms \cite{brevsar2013bucolic,genevois2017cubical}, and CNPC of spaces made of orthoschemes is closely related to swm-graphs \cite{chalopin_chepoi_hirai_osajda} and Helly graphs \cite{haettel2021lattices}. 
This raises two questions. 
First, whether there is a form of CNPC which is a common generalization of all these notions, hence can be applied to spaces made of mixed types of shapes. Second, whether each affine Coxeter group hints at a particular form of CNPC, which is applicable to complexes built with fundamental domains of such Coxeter group. This question is already unknown for Coxeter groups of type $\widetilde A_n$ with $n\ge 3$, which partially motivates this article. 

Attempting to answer the first question leads to the notion of \emph{weakly modular graph}. This notion, initially introduced in \cite{chepoi1989classification,bandelt1996helly}, demands metric balls in the graph satisfy a weak form of convexity (see Section~\ref{subsec:wm} for a detailed definition). It is a common generalization of bridged graphs, median graphs and Helly graphs mentioned in the previous paragraph, and serves as a mother notion to study its various sub-classes in \cite{chalopin_chepoi_hirai_osajda}. Though it is widely open whether weak modularity is compatible with other Coxeter shapes. Notable features of weakly modular graphs are that they enjoy a local-to-global characterization, and that they admit Euclidean isoperimetric inequalities.

Back to the $\widetilde A_n$ case of the second question with $n\ge 3$, as an initial step, it is shown by Munro \cite{munro2019weak} that while the 1-skeleton of the Coxeter complex of type $\widetilde A_n$ fails most form of CNPC, it does satisfy weak modularity. He also proved that the 1-skeleton of a 3-dimensional Euclidean building is weakly modular, though the high dimensional case remains open.
In this article we prove weak modularity for a much wider class of simplicial complexes with $\widetilde A_n$ simplices, including the higher dimensional $\widetilde A_n$ buildings. 

It turns out that many simplicial complexes may be endowed with natural $\widetilde A_n$ simplices, notably $\widetilde A_n$ Euclidean buildings, the Artin complex of the Artin-Tits group of type $\widetilde A_n$, and quotients of Garside groups. One common feature in all these examples is that the complex in question may be realized as the quotient of a (poset-theoretic) lattice under an action by $\Z$, as in \cite{haettel2021lattices} and \cite{haettel_injective_buildings}. The geometry of the corresponding lattice can be turned into a Helly graph by thickening, i.e. adding extra edges. However, there are no results about the original simplicial complex, i.e. the quotient of the lattice. 

For instance, Hoda \cite{hoda:crystallographic} proved that affine Coxeter groups of type $\widetilde A_n$ are not Helly, and Haettel \cite{haettel_injective_buildings} proved that if $n \geq 4$, then Euclidean buildings of type $\widetilde A_n$, even after equivariant thickening, are not Helly graphs.

In this article, we prove the following.

\bmthm (Theorem~\ref{thm:euclidean_building_wm})
The $1$-skeleton of any Euclidean building of type $\tilde{A}_n$ is a weakly modular graph.
\emthm

This theorem, as well as several other theorems below is a consequence of more general theorems on weak modularity of graphs coming from lattices with increasing $\Z$-action and the associated quotients, see Theorem~\ref{thm:quotient wm} and Theorem~\ref{theo:lattice wm}. We also note that the 1-skeletons considered in the above theorem also satisfy a stronger version of weak modular graph, as discussed in the end of Section~\ref{sec:quotient}.

On the other hand, the affine Coxeter complexes of type $\widetilde C_2$ and $\widetilde G_2$ are not weakly modular. Nevertheless, up to equivariantly adding edges, they become weakly modular. Concerning more general buildings, we formulate the following.

\bmconj \label{conj:building_wm}
Any building has an equivariant thickening which is a weakly modular graph.
\emconj

To be precise, we say that a graph $\Gamma$ is an equivariant thickening of a Euclidean building $X$ if $\Gamma$ contains $X$ as a subgraph, $X$ is quasi-isometric to $\Gamma$, and the automorphism group of $X$ extends as an automorphism group of $\Gamma$. This is motivated by the following particular cases.

\bmthm (Theorem~\ref{thm:many_buildings_thickening_wm})
\label{mainthm:many_buildings_thickening_wm}
Conjecture~\ref{conj:building_wm} holds for the following buildings:
\bit
\item Any spherical building.
\item Any Euclidean building of type $\widetilde A_n$, $\widetilde B_n$, $\widetilde C_n$, $\widetilde D_n$ or $\widetilde G_2$.
\item Any right-angled building.
\item Any rank $3$ building.
\item Any Gromov-hyperbolic building.
\eit
\emthm

We are also able to apply the same techniques to various classes of groups and complexes, first with Garside groups and weak Garside groups. Garside structures are essentially structures which locally look like a lattice, see Section~\ref{sec:garside} for precise definitions.

\bmthm (Theorem~\ref{thm:garside_cayley_wm})
Let $(G,\Delta,S)$ denote a (weak) Garside group. Then $\Cay(G,S)$, and its quotient by $\<\Delta\>$, are weakly modular graphs.
\emthm

This applies notably to finite type Garside groups, such as braid groups and spherical type Artin-Tits groups. Recall that Artin-Tits groups, defined by Tits in~\cite{tits_artin}, are natural generalizations of Coxeter groups and braid groups (see Section~\ref{subsec:Artin}). In this case of braid groups, one statement of the theorem is that the Cayley graph of braid groups with respect to simple braids is weakly modular. This also applies to infinite type Garside groups, such as braided crystallographic groups \cite{mccammond2017artin}, some Euclidean type Artin-Tits groups \cite{digne2006presentations,digne2012garside,mccammond2015dual} and some braid groups of imprimitive complex reflection groups \cite{corran2015braid}. Among weak Garside groups of finite types, one has all braid groups of complex reflection groups \cite{bessis2006non,bessis2015finite,corran2011new} except possibly the exceptional complex braid group of type $G_{31}$, all fundamental groups of complements of complexified real simplicial arrangements of hyperplane \cite{deligne1972immeubles},  and some extensions of Artin-Tits groups of type $B_n$ \cite{crisp2005representations}.

\mk

The techniques also apply to some Artin complexes. Recall that Artin-Tits groups have a natural candidate analogue of the curve complex known as the Artin complex, see~\cite{charney_davis_kpi1,cumplido2020parabolic}. It is the flag simplicial complex with vertices being cosets of maximal proper standard parabolic subgroups, with an edge for non-trivial intersection. In the case of Euclidean type Artin-Tits groups, the Artin complex is closely related to the Deligne complex.

\bmthm (Theorem~\ref{thm:Artin complex})
Let $A$ denote the Artin-Tits group of Euclidean type $\widetilde{A}_n$, and let $X$ denote the Artin complex of $A$. Then $X$ is a weakly modular graph.
\emthm

We observe the Artin complex of the Artin-Tits group of Euclidean type $\widetilde A_n$ has a topological interpretation as the complex of a certain collection of arcs in a surface (cf. Proposition~\ref{prop:iso}). Thus we give two treatments of Theorem~\ref{thm:Artin complex}, one is more on the Artin-Tits group side, the other uses surface topology and factors though the following theorem.

\bmthm (Theorem~\ref{thm:arc_complex_wm})
Let $n \geq 0$, and let $\Sigma$ be the $2$-sphere with $n+2$ punctures $N,S,p_1,\dots,p_n$. Let ${\cal A}(\Sigma)$ denote the subcomplex of the arc complex consisting of arcs between $N$ and $S$. Then ${\cal A}(\Sigma)$ is a weakly modular graph.
\emthm

The key objects in this article are lattices. In some applications, the corresponding lattice will be transparent: the lattice of norms in the case of $\widetilde A_n$ buildings, and the lattice of a Garside group. In some other examples, the lattice will be revealed through a mere local lattice property, as in the case of Artin complexes and arc complexes. To this end, we also reformulate work of Bessis into a simple local-to-global property for lattices, in the framework of Garside categories (see Section~\ref{subsec:lattice_local_global} for definitions).

\bmthm (Theorem~\ref{thm:lattice ltg})
Suppose $(P,\le)$ is a homogeneous weakly ordered set. If there exists an automorphism $\varphi: P\to P$ which generates $\le$ such that 
\begin{enumerate}
	\item $\varphi(x)>x$ for any $x\in P$;
	\item $X_\varphi$ is simply connected;
	\item $[x,\varphi(x)]$ is a lattice for any $x\in P$.
\end{enumerate}
Then $\le$ generates a partial order $\le_t$ on $P$, and $(P_{\ge_t x},\le_t)$ and $(P_{\le_t x},\le_t)$ are lattices for any $x\in P$.
\emthm

Interestingly, the framework of Garside categories in Bessis \cite{bessis2006garside} and the framework of lattices with $\Z$-action in Haettel \cite{haettel2021lattices} have many connections. Actually the former is a special case of the latter in an appropriate sense, and the latter contains a ``continuous'' version of the former (see the example and remark after Theorem~\ref{thm:building}). To this end, we record the following theorem, which gives a dictionary between different settings, under appropriate assumptions. See Section~\ref{subsec:dictionary} for precise definitions for terms in the following theorem. We hope this dictionary will help more geometric group theorists get interested in Garside groups.
\bmthm (Theorem~\ref{thm:dictionary_garside})
The following objects are equivalent:
\bit
\item A categorical Garside structure.
\item A Garside lattice.
\item A Garside flag complex.
\eit
\emthm

\mk

Note that, as most examples of Garside groups came from braid groups, the following questions were asked (see~\cite[Question~2]{garside}, \cite{calvez_wiest_1,calvez_wiest_2}):

\bq
Are Garside groups linear? Do Garside groups act on Gromov-hyperbolic spaces?
\eq

We use our framework to promote lattices in triangular buildings into Garside groups with exotic properties. This is one of the first and simplest evidence that affine buildings are closely related to Garside theory. In particular, we answer the previous questions by the negative, in the following way.

\bmthm (Corollary~\ref{cor:exotic_garside})\

There exists a non-linear Garside group $G_1$.

There exists a Garside group $G_2$ without any non-trivial action on a Gromov-hyperbolic space. Moreover, any normal subgroup of $G_2$ is either finite, or has virtually cyclic quotient.

In addition, for each of these Garside groups $G_i$, the central quotient $G_i/Z(G_i)$ has Lafforgue's strong Property (T).
\emthm

Note that, for the Garside group $G_2$, the hyperbolic space constructed by Calvez and Wiest in~\cite{calvez_wiest_1}, called the additional length complex, has finite diameter.

\mk

{\bf Acknowledgements:} The authors thank warmly to the anomalous referee for a very careful reading of the article and many helpful comments.
The authors would like to thank Andrew Putman for precise references concerning the arc complex on the punctured sphere. The authors would also like to thank Owen Garnier and Ivan Marin for many precisions concerning complex braid groups. The authors thank Anthony Genevois for many helpful comments. The authors thank Damian Osajda for stimulating discussions. The authors thank Sami Douba and Jean L\'ecureux for discussions on buildings.

\tableofcontents

\section{Preliminary}

\subsection{Artin-Tits groups and Coxeter groups}
\label{subsec:Artin}
Let $\Gamma$ be a finite simple graph with each edge labeled by an integer $\ge 2$. The \emph{Artin-Tits group with defining graph $\Gamma$}, also known as Artin group, denoted $A_\Gamma$, is given by the following presentation. Generators of $A_\Gamma$ are in one to one correspondence with vertices of $\Gamma$, and there is a relation of the form
\begin{center}
	$\underbrace{aba\cdots}_{m}=\underbrace{bab\cdots}_{m}$
\end{center}
whenever two vertices $a$ and $b$ are connected by an edge labeled by $m$. The \emph{Coxeter group with defining graph $\Gamma$}, denoted $W_\Gamma$, has the same generator sets and the same relators as the Artin-Tits group, with extra relations $v^2=1$ for each vertex $v\in\Gamma$.

An Artin group is \emph{spherical}, if the associated Coxeter group is finite. 
A \emph{standard parabolic subgroup} of $A_\Gamma$ is a subgroup generated by a subset of the vertices of $\Gamma$. Each standard parabolic subgroup is an Artin group as well \cite{vanderlek}. A \emph{parabolic subgroup} is a conjugate of a standard parabolic subgroup.

In this article we are mostly interested in Artin groups and Coxeter groups of type $\widetilde A_{n-1}$, in which case the defining graph $\Gamma$ is a complete graph on $n$ vertices such that there exists an embedded $n$-cycle in $\Gamma$ with each edge contained in the cycle are labeled by $3$ and all other edges are labeled by $2$.

Let $A_{\Gamma}$ be an Artin-Tits group of type $\widetilde A_{n-1}$ with consecutive generators in its Dynkin diagram labeled by $\{s_1,s_2,\ldots, s_n\}$. Let $A_i$ be the subgroup of $A_\Gamma$ generated by all the generators except $s_i$. Let $X$ be the \emph{Artin complex} of $A_\Gamma$ (cf. \cite{charney_davis_kpi1,cumplido2020parabolic}), namely vertices of $X$ are in 1-1 correspondence with left cosets of the form $\{gA_i\}_{g\in A_\Gamma, 1\le i\le n}$. Two vertices are joined by an edge if the associated left cosets have non-empty intersection. Then $X$ is defined to be the flag completion of its 1-skeleton. 

Note that the barycentric subdivision of $X$ coincides with another complex, called the modified Deligne complex of $A_\Gamma$ as defined in \cite{charney_davis_kpi1}, since vertices of the modified Deligne complex correspond to left cosets of spherical standard parabolic subgroups, and in type $\widetilde A_{n-1}$ the spherical standard parabolic subgroups are exactly the proper standard parabolic subgroups.

\subsection{Weakly modular graphs and local to global}
\label{subsec:wm}
Let $\Gamma$ be a simplicial graph. We  endow $\Gamma$ with the path metric such that each edge has length $=1$.
Recall that $\Gamma$ is \emph{weakly modular} if for every vertex $x\in\Gamma$, and every positive integer $n\ge 2$, the following two conditions hold:
\begin{enumerate}
	\item (\emph{triangle condition} TC$(x)$) for any two adjacent vertices $x_1,x_2\in\Gamma$ such that $d(x,x_1)=d(x,x_2)=n$, there exists vertex $y$ with $d(y,x)=n-1$ such that $d(y,x_1)=d(y,x_2)=1$;
	\item (\emph{quadrangle condition} QC$(x)$) for for any two vertices $x_1,x_2\in\Gamma$ with $d(x,x_1)=d(x,x_2)=n$ and $d(x_1,x_2)=2$ such that $x_1$ and $x_2$ are adjacent to a common vertex at distance $n+1$ from $x$, there exists a vertex $y$ with $d(y,x)=n-1$ such that $d(y,x_1)=d(y,x_2)=1$.
\end{enumerate}

A graph $\Gamma$ is \emph{local weakly modular} if for every vertex $x\in\Gamma$, the triangle condition TC$(x)$ and quadrangle condition TC$(x)$ hold with $n=2$. 

An \emph{induced} subgraph of $\Gamma$ is a subgraph which contains all edges of $\Gamma$ that join two vertices of the subgraph. A \emph{square} of $\Gamma$ is an induced embedded 4-cycle of $\Gamma$.

Given a graph $\Gamma$, the \emph{triangle-square complex} of $\Gamma$ is a two-dimensional cell complex with 1-skeleton $\Gamma$ such that we fill in a solid triangle for each embedded 3-cycle in the graph and fill in a solid square for each square of $\Gamma$.

\bthm(\cite[Theorem 3.1]{chalopin_chepoi_hirai_osajda})
\label{theo:wm ltg}
If $\Gamma$ is a local weakly modular graph, and its triangle-square complex is simply-connected, then $\Gamma$ is a weakly modular graph.
\ethm

Actually the converse to this theorem is also true, i.e. if $\Gamma$ is a weakly modular graph, then its triangle-square complex is simply-connected \cite[Lemma 5.5]{brevsar2013bucolic}.
\subsection{Posets, Lattices and local to global}

\label{subsec:lattice_local_global}

Let $P$ be a poset, i.e. a partially ordered set.
Let $S\subset P$. An \emph{upper bound} (resp. lower bound) for $S$ is an element $x\in P$ such that $s\le x$ (resp. $s\ge x$) for any $s\in S$. The \emph{join} of $S$ is an upper bound $x$ of $S$ such that $x\le y$ for any other upper bound $y$ of $S$. The \emph{meet} of $S$ is a lower bound  $x$ of $S$ such that $x\ge y$ for any other lower bound $y$ of $S$. We will write $x\vee y$ for the join of two elements $x$ and $y$, and $x\wedge y$ for the meet of two elements (if the join or the meet exists).
We say $P$ is \emph{lattice} if $P$ is a poset and any two elements in $P$ have a join and have a meet. For $a,b\in P$ with $a\le b$, the \emph{interval}  between $a$ and $b$, denoted by $[a,b]$, is the collection of all elements $x$ of $P$ such that $a\le x$ and $x\le b$. A poset $P$ is homogeneous if there is a function $\ell$ from each comparable pair in $P$ to the non-negative integers such that if $a\le b\le c$, then $\ell(a\le c)=\ell(a\le b)+\ell(b\le c)$; and $\ell(a\le b)=0$ if and only if $a=b$. Note that if $P$ is homogeneous lattice, then any upper bounded subset of $P$ has a join and any lower bounded subset of $P$ has a meet. We will also need the following notion of \emph{weakly ordered} set which generalizes the notion of poset by allowing the transitivity to fail.
\begin{defi}
	\label{defi:weak_order}
A \emph{weakly ordered} set $P$ is a set with a binary relation $\le$ over $P$ which is  reflexive and antisymmetric. Moreover, while transitivity may fail, we do require the following associativity law for transitivity. Define $(a,b,c)\in P^3$ to be a \emph{transitive triple} if $a\le b$, $b\le c$ and $a\le c$. We require $\le$ satisfies the following condition:

\noindent
$(\ast)$: for any quadruple $a,b,c,d\in P$ with $a\le b,b\le c$ and $c\le d$, we have $(b,c,d)$ and $(a,b,d)$ are transitive triples if and only if $(a,b,c)$ and $(a,c,d)$ are transitive triples.
\end{defi}

The notions of upper bound, join, lower bound, meet, interval and homogeneity can be defined for a weakly ordered set in the same way. 
Let $(P,\le)$ be a weakly ordered set. For $x\in P$, let $P_{\ge x}$ be the collection of all elements which are $\ge x$. Similarly we define $P_{\le x}$. Note that $P_{\ge x}$ and $P_{\le x}$ are actually posets by $(\ast)$. 
A \emph{weak chain} in a weakly ordered set $P$ is a sequence of elements $c_1<c_2<c_3<\cdots<c_n$ such that any two adjacent elements in the sequence are comparable. A \emph{chain} is a weak chain such that any two elements in the weak chain are comparable.
Note that if $P$ does not contain non-trivial weak chains which start and end at the same element, then the weak order $\le$ on $P$ actually generates a partial order $\le_t$, where $a\le_tb$ if $a$ and $b$ are the first and the last member of a weak chain in $P$ with respect to $\le$.

An \emph{automorphism} of $(P,\le)$ is a bijection of $P$ preserving the relation $\le$.
Suppose $(P,\le)$ is a weakly ordered set with an automorphism $\varphi: P\to P$ such that $\varphi(x)>x$ for any $x\in P$. Let $\le_{\varphi}$ be the subrelation of $\le$ made of all possible $a\le b$ such that there exists $x\in P$ such that $a,b\in[x,\varphi(x)]$. We say $\varphi$ \emph{generates} $\le$ if $a$ and $b$ fit into the first element and the last element of a weak chain with respect to $\le_{\varphi}$ whenever $a\le b$.

Let $X_\varphi$ be the simplicial complex whose vertex set is $P$, and $a$ and $b$ are joined by an edge if $a\le_\varphi b$. Simplices of $X_\varphi$ correspond to weak chains in $P$ that are contained in an interval of the form $[x,\varphi(x)]$ for some $x\in P$. Note that though a weak chain contained in $[x,\varphi(x)]$ is automatically a chain.

The following is a consequence of work of Bessis \cite{bessis2006garside}.
\bthm
\label{thm:lattice ltg}
Suppose $(P,\le)$ is a homogeneous weakly ordered set. If there exists an automorphism $\varphi: P\to P$ which generates $\le$ such that 
\begin{enumerate}
	\item $\varphi(x)>x$ for any $x\in P$;
	\item $X_\varphi$ is simply connected;
	\item $[x,\varphi(x)]$ is a lattice for any $x\in P$.
\end{enumerate}
Then $\le$ generates a partial order $\le_t$ on $P$, and $(P_{\ge_t x},\le_t)$ and $(P_{\le_t x},\le_t)$ are homogeneous lattices for any $x\in P$. Moreover, $X_\varphi$ is contractible.
\ethm

\begin{proof}
By definition of $\le_\varphi$, we have $x\le a\le \varphi(x)$ if and only if $x\le_\varphi a\le_\varphi \varphi(x)$. Thus the interval $[x,\varphi(x)]$ with respect to $\le$ and the same interval with respect to $\le_\varphi$ are the same. 
Thus we can assume without loss of generality that $\le=\le_\varphi$. We claim for any $x\in P$, if $x\le a$, then $a\le \varphi(x)$. Indeed, $x\le a$ implies there exists $y\in P$ such that $x,a\in [y,\varphi(y)]$. As $y\le x$, we have $\varphi(y)\le\varphi(x)$. We now consider the quadruple $x,a,\varphi(y),\varphi(x)$. Then $(x,a,\varphi(y))$ and $(x,\varphi(y),\varphi(x))$ are transitive triples as $x\in [y,\varphi(y)]$. Thus $a\le\varphi(x)$ and the claim is proved. Similarly, we know for any $x\in P$, if $a\le x$, then $\varphi^{-1}(x)\le a$.

Note that $(P,\le)$ is a special example of a germ in the sense of \cite[Definition 1.1]{bessis2006garside}. More precisely, the objects of the germ are elements in $P$, and there is a morphism from $a$ to $b$ if $a\le b$. This germ is homogeneous Garside in the sense of \cite[Definition 3.2]{bessis2006garside}. Indeed, \cite[Definition 3.2]{bessis2006garside} (i) is clear. By the above claim, the collection of morphism starting at $x$ has a maximal element, which is $x\le \varphi(x)$, thus \cite[Definition 3.2]{bessis2006garside} (ii) holds true. Now \cite[Definition 3.2]{bessis2006garside} (iv) follows from assumption (3). \cite[Definition 3.2]{bessis2006garside} (iii) translates into the following: we consider the map from the set of morphisms starting at $x$ to the set of morphisms ending at $\varphi(x)$ sending $x\le a$ to $a\le \varphi(x)$. By the above claim, this map is well-defined and is a bijection. Now \cite[Definition 3.2]{bessis2006garside} (iii) is clear.

Consider the category $\mathcal C$ generated by $\le$, whose objects are $P$ and whose morphisms are equivalent classes of chains of the form $a_1\le a_2\le \cdots \le a_n$ (we require adjacent members in the sequence to be comparable, however, non-adjacent members might not be comparable), where the equivalence relation is generated by $\sim$ with $$(a_1\le \cdots a_i\le a_{i+1}\le a_{i+2} \le \cdots \le a_n)\sim (a_1\le \cdots a_i\le a_{i+2} \le \cdots \le a_n)$$ if either $(a_i,a_{i+1},a_{i+2})$ is a transitive triple or $a_i=a_{i+1}$.

By \cite[Theorem 3.3]{bessis2006garside}, $\mathcal C$ is a categorical Garside structure in the sense of \cite[Definition 2.4]{bessis2006garside}. In particular, the collection of all morphisms starting from $x\in P$, endowed with the prefix order (i.e. $f\le_p g$ if $g=fh$ for a morphism $h$ of $\mathcal C$), is a lattice; and the collection of all morphisms ending at $x\in P$, endowed with the suffix order (i.e. $f\ge_s g$ if $hg=f$ for a morphism $h$ of $\mathcal C$), is a lattice. These
lattices are homogeneous by \cite[Theorem 3.3]{bessis2006garside}, where the length function $\ell(f\le_p g)$ is defined to be the sum $\sum_{i=1}^n \ell(f_i\le_p f_{i+1})$ where $f_1\le_p f_2\le_p\cdots\le_p f_n$ is a chain with $f_1=f$ and $f_n=g$ such that adjacent members in the chain are comparable in $P$.
Thus we are done as long as we can show there are no non-trivial weak chains in $P$ which start and end at the same element.

Let $\mathcal G$ be the groupoid obtained by adding formal inverses to all morphisms in $\mathcal C$. It follows from \cite[Section 2]{bessis2006garside} (more precisely the discussion of normal forms of morphisms in $\mathcal G$ between \cite[Definition 2.10]{bessis2006garside} and \cite[Definition 2.11]{bessis2006garside}) that the map $\mathcal C\to \mathcal G$ is injective.
Take an object $x\in P$, let $\mathcal G_x$ be the collection of morphisms in $\mathcal G$ starting at $x$. Let $|\mathcal G_x|$ be the simplicial complex defined by Bessis as follows. The vertices are corresponding to elements in $|\mathcal G_x|$. Two different vertices are joined by an edge if the corresponding two morphisms $f$ and $g$ satisfy $f=gh$ or $g=fh$ for some simple morphism $h$ in $\mathcal C$
(a \emph{simple} morphism of $\mathcal C$ is a morphism of the form $a<b$ such that $b\le \varphi(a)$). Then $|\mathcal G_x|$ is the flag complex of its 1-skeleton. By \cite[Corollary 7.6]{bessis2006garside}, $|\mathcal G_x|$ is contractible, hence simply-connected. Note that there is a map $p$ from $\mathcal G_x$ to $P$ by sending each morphism to its endpoint. By the proof of \cite[Corollary 7.6]{bessis2006garside}, a simplex in $|\mathcal G_x|$ corresponds to a collection $\{f_1,\ldots,f_k\}$ such that $f_j$ is the composition of $f_i$ with some simple morphism in $\mathcal C$ whenever $j>i$. Thus the map $p$ extends to a simplicial map $p:|\mathcal G_x|\to X_\varphi$ which is also a covering map. As $X_\varphi$ is simply connected, we know $p$ is a simplicial isomorphism. If there exist a non-trivial weak chain in $P$ starting and ending at the same element, then there exists a morphism $g$ in $\mathcal G_x$, and a collection of non-trivial simple morphisms $\{h_1,\ldots,h_k\}$ in $\mathcal C$ such that $gh_1\cdots h_k=g$. By cancellation property in the groupoid $\mathcal G$, $h_1\cdots h_k$ is an identity morphism in $\mathcal G$, hence in $\mathcal C$ due to the injectivity of $\mathcal C\to \mathcal G$. However, this contradicts the homogeneity assumption on $\mathcal C$.
\end{proof}

\section{The diagonal quotient of a lattice is weakly modular}
\label{sec:quotient}
Assume that $L$ is a lattice, such that each nonempty upper bounded subset of $L$ has a join (as a consequence, each nonempty lower bounded subset of $L$ has a meet). 
Assume that there is an order-preserving increasing action of $\Z$ on $L$, notated additively (i.e. the image of $x\in L$ under the action $n\in\Z$ is denoted by $x+n$), such that
$$\forall x,y \in L, \exists k \in \N, x-k \leq y \leq x+k.$$

\mk

We will define a graph $X$ from $L$, with vertex set $L/\Z$. Add an edge between $x,y \in X$ if, for some representatives $x_0,y_0$ of $x,y$ in $L$, we have
$$ x_0 \leq y_0 \leq x_0+1.$$

We have the following.

\bthm
\label{thm:quotient wm}
$X$ is a weakly modular graph.
\ethm

Note that $X$ is connected by Lemma~\ref{lem:distance} below and our assumption on $L$. Thus the theorem is a consequence of Lemma~\ref{lem:triangle} and Lemma~\ref{lem:quadrangle} below.

\blem
\label{lem:distance}
Given any vertices $x,y \in X$ and any representatives $x_0,y_0 \in L$, we have
$$d(x,y) = \min\{n \geq 0 \st \exists k,h \in \Z, x_0+k \leq y_0+h \leq x_0+k+n\}.$$
\elem

\bp
Let us denote the formula in the statement by $d'$. Given an edge path in $X$, we can consider a lift of consecutive vertices in this path to $L$ of the form $x_0 \leq x_1 \leq \dots \leq x_n$ such that $x_i\le x_{i+1}\le x_i+1$ for all $i$. Then $x_0 \leq x_n \leq x_0+n$, so $d'(x_0+\Z,x_n+\Z) \leq n$. Hence $d' \leq d$.

\mk

Conversely, we will prove by induction on $n \geq 0$ that, if $x,y \in L$ are such that $d'(x,y)=n$, then $d(x,y)=n$. When $n \leq 1$ it is obvious. Assume that the statement holds for values smaller than $n$, and fix vertices $x,y \in X$ and representatives $x_0,y_0 \in L$ such that $x_0 \leq y_0 \leq x_0+n$, where $n=d'(x,y)$.

\mk

Let $z_0 = (x_0+1) \wedge y_0$: we have $z_0+n-1 =(x_0+n) \wedge (y_0+n-1)$. As $y_0 \leq x_0+n$ and $y_0 \leq y_0+n-1$, we have $z_0 \leq y_0 \leq z_0+n-1$, so $d'(y,z) \leq n-1$. By induction, we deduce that $d(y,z)=d'(y,z) \leq n-1$. Furthermore, we have $d(x,z) \leq 1$, so as $d(x,y) \geq d'(x,y)=n$ we conclude that $d(x,y)=n=d'(x,y)$.
\ep

\blem
\label{lem:triangle}
$X$ satisfies the triangle condition: fix $n \geq 2$, and let $x,y,z \in X$ such that $d(x,y)=d(x,z)=n$ and $d(y,z)=1$. There exists $u \in X$ such that $d(x,u)=d(u,y)=1$ and $d(u,z)=n-1$.
\elem

\bp
Consider representatives $x_0,y_0,z_0$ such that $x_0 \leq y_0 \leq x_0+n$ and $y_0 \leq z_0 \leq y_0+1$. We will prove that we can furthermore assume that $z_0 \leq x_0+n$.

\mk

Since $x_0 \leq z_0 \leq y_0+1 \leq x_0+n+1$ and $d(x,z)=n$, by Lemma~\ref{lem:distance} there are two only possibilities (note that ):
\bit
\item either $x_0 \leq z_0 \leq x_0+n$
\item or $x_0+1 \leq z_0 \leq x_0+n+1$.
\eit
To see these are the only two possibilities, note that $d(x,z)=n$ implies that $[x_0,x_0+n]$ contains a unique lift of $z$.

In the former case, we have indeed $z_0 \leq x_0+n$. In the latter case, we have $x_0 \leq z_0-1 \leq y_0 \leq x_0+n$ and $z_0-1 \leq y_0 \leq z_0$, so up to replacing the pair $(y_0,z_0)$ by the pair $(z_0-1,y_0)$, we can always assume that $z_0 \leq x_0+n$.

\mk

Let $u_0=(x_0 +n-1) \wedge y_0$: we have $x_0 \leq u_0 \leq x_0+n-1$, so $d(x,u) \leq n-1$.

\mk

Furthermore, since $z_0 \leq x_0+n$ and $z_0 \leq y_0+1$, we deduce that $u_0+1 = (x_0+n) \wedge (y_0+1) \geq z_0$, hence $u_0 \leq y_0 \leq z_0 \leq u_0+1$. Hence $d(u,y) \leq 1$ and $d(u,z) \leq 1$.

\mk

Hence we have $d(x,u)=n-1$ and $d(u,y)=d(u,z)=1$.
\ep

\blem
\label{lem:quadrangle}
$X$ satisfies the quadrangle condition: let $n \geq 2$ and $x,y,z,t \in X$ such that $d(x,y)=d(x,z)=n$, $d(y,t)=d(z,t)=1$ and $d(x,t)=n+1$. There exists $u \in X$ such that $d(x,u)=n-1$ and $d(u,y)=d(u,z)=1$.
\elem

\bp
We first prove that there exists $s \in X$ such that $d(y,s)=d(z,s)=1$ and $d(x,s) \leq n$. Consider representatives $x_0,y_0,z_0,t_0$ such that $x_0 \leq y_0 \leq x_0+n$, $x_0 \leq z_0 \leq x_0+n$ and $y_0 \leq t_0 \leq y_0+1$.

\mk

We will prove that also $z_0 \leq t_0 \leq z_0+1$. Since $x_0 \leq y_0 \leq t_0 \leq y_0+1 \leq x_0+n+1$ and $d(x,t)=n+1$, we deduce from Lemma~\ref{lem:distance} that $t_0$ is not comparable to $x_0+1$ nor $x_0+n$. Since $z_0 \leq x_0+n$, we deduce that $t_0$ is not inferior to $z_0$, hence $t_0 \geq z_0$ as $d(z,t)=1$. Similarly, since $z_0+1 \geq x_0+1$, we deduce that $t_0 \leq z_0+1$. So we have $z_0 \leq t_0 \leq z_0+1$.

\mk

Let $s_0=y_0 \wedge z_0$. Since $x_0 \leq s_0 \leq x_0+n$, we know that $d(x,s) \leq n$. Furthermore, since $t_0-1 \leq y_0\leq t_0$ and $t_0-1\leq z_0 \leq t_0$, we deduce that also $t_0-1 \leq s_0 \leq t_0$. Thus if we take $s\in X$ to be the vertex associated with $s_0$, then $d(s,y) \leq 1$ and $d(s,z) \leq 1$.

\mk

Let $u_0=(x_0+n-1) \wedge s_0$: we have $d(x,u) \leq n-1$. Moreover, $u_0\le y_0$ and $u_0\le z_0$. We will prove that $y_0\leq u_0+1$ and $z_0 \leq u_0+1$. Since $s_0=y_0 \wedge z_0$, we have $u_0 = (x_0+n-1) \wedge y_0 \wedge z_0$, so $u_0+1 = (x_0+n) \wedge (y_0+1) \wedge (z_0+1)$.
Note that $y_0 \leq x_0+n$ and $y_0 \leq y_0+1$. Furthermore, we have $y_0 \leq t_0 \leq z_0+1$. As a conclusion, we have $y_0 \leq u_0+1$, and similarly $z_0 \leq u_0+1$. We conclude that $d(u,y) \leq 1$ and $d(u,z) \leq 1$.

\mk

Hence $d(x,u)=n-1$ and $d(u,y)=d(u,z)=1$.
\ep

\begin{rmk}
The above argument implies that $X$ satisfies the following  stronger versions of the triangle condition and the quadrangle condition, namely:
\begin{enumerate}
	\item for any vertex $x\in X$ and any complete subgraph $Y\subset X$ such that each vertex of $Y$ is at distance $n$ from $x$, there exists a vertex $z\in X$ such that $d(x,z)=n-1$ and $z$ is adjacent to each vertex in $Y$;
	\item for any vertex $x\in X$ and any vertex $t\in X$ with $d(x,t)=n+1$, let $Y=\{y\in X:d(y,t)=1\ \textrm{and}\ d(y,x)=n\}$, then there exists a vertex $u$ at distance $n-1$ from $x$ such that $u$ is adjacent to each vertex in $Y$, moreover, there exists a vertex $s\in Y$ such that $s$ is adjacent to each vertex in $Y\setminus\{s\}$.
\end{enumerate}
The proof of property (2) is identical to Lemma~\ref{lem:quadrangle}. The proof of property (1) is a small adjustment of Lemma~\ref{lem:triangle}. Namely take $y\in Y$ and choose representatives $x_0,y_0$ of $x,y$ such that $x_0\le y_0\le x_0+n$. Then for any $z\in Y$, the proof of Lemma~\ref{lem:triangle} implies that there exists a representative $z_0$ of $z$ such that either $y_0-1\le z_0\le y_0$ and $x_0\le z_0\le x_0+n$ or $y_0\le z_0\le y_0+1$ and $x_0\le z_0\le x_0+n$. 
Let $Y_0$ be the collection of all such representatives of elements in $Y$. We claim each pair of elements $z_0,z'_0$ in $Y_0$ are comparable, and if $z_0\ge z'_0$, then $z_0\le z'_0+1$. First we consider the case $z_0\ge y_0\ge z'_0$. The first part of the claim is clear. As $z$ and $z'$ are adjacent, $z_0$ and $z'_0+1$ are comparable. If $z_0>z'_0+1$, then $x_0+n-1\ge z_0-1\ge z'_0\ge x_0$, which contradicts that $d(z,x)=n$. Thus $z_0\le z'_0+1$. Now we consider the case where both $z_0,z'_0$ are lower bounded by $y_0$, then $z_0,z'_0$ are upper bounded by $y_0+1$. This, together with the fact that $z$ and $z'$ are adjacent imply the claim. The case $z_0,z'_0$ being lower bounded by $y_0$ is similar. Thus the claim is proved.
Let $v_0$ be the meet of $Y_0$ and $u_0=(x_0+n-1)\wedge v_0$. By the claim, for any $z_0\in Y_0$, we have $z_0-1\le v_0\le z_0$. In particular, $z_0\le v_0+1$. The argument in Lemma~\ref{lem:triangle} implies that $u_0\le v_0\le z_0\le u_0+1$. Hence $d(u,z)\le 1$ for any $z\in Y$. 
\end{rmk}

\section{A lattice with a diagonal action is weakly modular}
\label{sec:lattice_no_quotient}

Assume that $L$ is a lattice, such that each nonempty upper bounded subset of $L$ has a join. Assume that there is an order-preserving increasing action of $\Z$ on $L$, notated additively, such that
$$\forall x,y \in L, \exists k \in \N, x-k \leq y \leq x+k.$$

\mk

We will define a graph $X$ from $L$, with vertex set $L$. Add an edge between $x,y \in X$ if $x \leq y \leq x+1$ or $y \leq x \leq y+1$.

\bthm
\label{theo:lattice wm}
$X$ is a weakly modular graph.
\ethm

This theorem is a consequence of Theorem~\ref{theo:wm ltg}, as well as Lemma~\ref{lem:connected}, Lemma~\ref{lem:tri}, Lemma~\ref{lem:quad} and Lemma~\ref{lem:sc} below.

\blem
\label{lem:connected}
The graph $X$ is connected.
\elem

\bp
Fix $x,y \in X$. Since there is an edge between $x$ and $x+1$, we may assume that $x \leq y$. For each $n \in \N$, let $x_n=(x+n) \wedge y$. For each $n \in \N$, we have $x_n = (x+n) \wedge y \leq (x+n+1) \wedge y = x_{n+1}$, and also $x_{n+1} = (x+n+1) \wedge y \leq (x+n+1) \wedge (y+1) = ((x+n) \wedge y)+1=x_n+1$. So $x_n$ and $x_{n+1}$ are adjacent in $X$.

By assumption, there exists $n \in \N$ such that $y \leq x+n$. So we deduce that $x_n=y$, and $x$ and $y$ are connected by the path $x_0=x,x_1,\dots,x_n=y$ in $X$.
\ep

\blem
Given any vertices $x,y \in X$, we have
$$d(x,y) = \min\{n+m \st n,m \in \N, x \leq y+n, y \leq x+m\}.$$
Moreover, $x \vee y$ and $x \wedge y$ each belong to a geodesic between $x$ and $y$.
\elem

\bp
Fix $x,y \in X$, and consider minimal $n,m \in \N$ such that $x \leq y+n, y \leq x+m$.

\mk

According to the previous proof, there is a length $n$ path from $x \wedge y$ to $x \wedge (y+n)=x$. There is also a length $m$ path from $x \wedge y$ to $(x+m) \wedge y=y$. Hence $d(x,y) \leq n+m$.

\mk

Conversely, we will prove by induction on $n+m$ that $d(x,y)=n+m$. Let $y' \in X$ be adjacent to $y$, we will prove that the corresponding integers for $x$ and $y'$ satisfy $n'+m' \leq n+m+1$.

If $y \leq y' \leq y+1$, then $x \leq y+n \leq y'+n$ and $y' \leq y+1 \leq x+m+1$, so $n' \leq n$ and $m' \leq m+1$.

If $y' \leq y \leq y'+1$, then $x \leq y+n \leq y'+n+1$ and $y' \leq y \leq x+m$, so $n' \leq n+1$ and $m' \leq m$.

\mk

We conclude that $d(x,y)=n+m$.
\ep

\blem
\label{lem:tri}
$X$ satisfies the triangle condition: fix $n \geq 2$, and let $x,y,z \in X$ such that $d(x,y)=d(x,z)=n$ and $d(y,z)=1$. There exists $u \in X$ such that $d(u,y)=d(u,z)=1$ and $d(u,x)=n-1$.
\elem

\bp
Assume, without loss of generality, that $y \leq z \leq y+1$. Let $p,m \in \N$ be minimal such that $x \leq y+p$ and $y \leq x+m$. 

Since $d(x,z)=d(x,y)=p+m$, there are two possibilities:
\bit
\item either $x \leq z+p$ and $z \leq x+m$. Assume first that $m \geq 1$, and let $u = y \wedge (x+m-1)$: we have $d(u,y \wedge x) \leq m-1$ and $d(y \wedge x,x) \leq p$, so $d(u,x) \leq p+m-1$. Also $u \leq y \leq z$, and $y,z \leq y+1, x+m$, so $y,z \leq u+1$: we have $d(u,y) \leq 1$ and $d(u,z) \leq 1$. By the triangular inequality we have $d(x,u)=p+m-1$ and $d(u,y)=d(u,z)=1$.

Assume now that $m=0$, let us then define $u = (y+1) \wedge x$. We have $u \leq x$ and $x \leq y+p \leq u + p-1$, hence $d(u,x) \leq p-1$. Also $u \geq z \wedge x=z \geq y$, and $u \leq y+1 \leq z+1$, so $d(u,y) \leq 1$ and $d(u,z) \leq 1$. By the triangular inequality we have $d(x,u)=m-1$ and $d(u,y)=d(u,z)=1$.
\item or $x \leq z+p-1$ and $z \leq x+m+1$.

Let $u = z \wedge (x+m)$: we have $d(u,z \wedge x) \leq m$ and $d(z \wedge x,x) \leq p-1$, so $d(u,x) \leq p+m-1$. Also $y \leq z,x+m$ so $y \leq u$. And $u \leq z \leq y+1$, so $d(u,y) \leq 1$. We also have $u \leq z$ and $z \leq (z+1), (x+m+1)$ so $z \leq  u+1$: we have $d(u,z) \leq 1$. By the triangular inequality we have $d(x,u)=p+m-1$ and $d(u,y)=d(u,z)=1$.
\eit
\ep

\blem
\label{lem:quad}
$X$ satisfies the local quadrangle condition. More precisely, for any $x,y,z,t \in X$ such that $d(x,y)=d(x,z)=2$, $d(x,t)=3$ and $d(y,t)=d(z,t)=1$, there exists $u \in X$ such that $d(x,u)=d(y,u)=d(z,u)=1$.
\elem

\bp
Note that, if $d(x,y)=2$, there are three possibilities: $x \leq y+1$ and $y \leq x+1$, $x \leq y \leq x+2$ and $y \leq x \leq y+2$. In this proof, we will call the last two possibilities of type $(2,0)$.

\bit
\item Let us first assume that $x \leq y+1$, $y \leq x+1$, $x \leq z+1$, $z \leq x+1$, $y \leq z+1$ and $z \leq y+1$. Let $u=x \wedge y \wedge z$. We have $u \leq x$ and, since $x \leq x+1, y+1, z+1$, we have $x \leq u+1$. So $d(x,u) \leq 1$. Similarly $d(y,u) \leq 1$ and $d(z,u) \leq 1$. By the triangular inequality we have $d(x,u)=d(y,u)=d(z,u)=1$.
\item Assume now that $x \leq y+1$, $y \leq x+1$, $x \leq z+1$, $z \leq x+1$ and $y \leq z \leq z+2$. We will show that this contradicts the existence of $t$. Since $d(x,t)=3$, there are three possibilities:
\ben
\item If $x \leq t \leq x+3$, then since $d(y,t)=d(z,t)=1$ we deduce that $t \leq z+1 \leq x+2$, which is a contradiction.
\item If $x-1 \leq t \leq x+2$, then since $z \leq x+1$, we have $t \not\leq z$, so $y \leq z \leq t$. Hence $t \leq y+1$, so $z \leq y+1$, which is a contradiction.
\item If $x-2 \leq t \leq x+1$, then since $x-1 \leq y$, we have $t \not\geq y$, so $t \leq y \leq z$. Hence $t \geq z-1$, so $z \leq y+1$, which is a contradiction.
\item If $x-3 \leq t \leq x$, then since $d(y,t)=d(z,t)=1$ we deduce that $t \geq z-1 \leq x-2$, which is a contradiction.
\een
\item Assume now that $x \leq y \leq x+2$, $x \leq z+1$, $z \leq x+1$, $y \leq z+1$ and $z \leq y+1$. We will show that this contradicts the existence of $t$. We know that $t \leq z+1 \leq x+2$ and $t \geq y-1 \geq x-1$. Hence $x-1 \leq t \leq x+2$. So $t \not\leq x+1$, hence $t \not\leq z$: $z \leq t$. Similarly $t \not\geq x$, hence $t \not\geq y$: $t \leq y$. We deduce that $z \leq t \leq y$, which contradicts $d(y,z)=2$.
\item Assume now that $y \leq x \leq y+2$, $x \leq z+1$, $z \leq x+1$, $y \leq z+1$ and $z \leq y+1$. We will show that this contradicts the existence of $t$. We know that $t \leq y+1 \leq x+1$ and $t \geq z-1 \geq x-2$. Hence $x-2 \leq t \leq x+1$. So $t \not\leq x$, hence $t \not\leq y$: we have $t \geq y$. Similarly $t \not\geq x-1$, hence $t \not\geq z$: we have $t \leq z$. We deduce $y \leq t \leq z$, which contradicts $d(y,z)=2$.
\item Assume now that there are two distances of type $(2,0)$: we will show the existence of $u$ independently of the assumption on $t$. Remark that if $x \leq y \leq x+2$ and $z \leq x \leq z+2$, we have $z \leq y \leq z+1$, which contradicts $d(y,z)=2$. So we may assume that, for instance, we have $x \leq y \leq x+2$, $x \leq z \leq x+2$, $y \leq z+1$ and $z \leq y+1$. Let $u=(x+1) \wedge y \wedge z$. Then $x \leq u \leq x+1$, so $d(u,x) \leq 1$. Also $u \leq y \leq u+1$, so $d(u,y) \leq 1$. Similarly $d(u,z) \leq 1$. By the triangular inequality we have $d(x,u)=d(y,u)=d(z,u)=1$.
\item Assume now that the three distances are of type $(2,0)$. Without loss of generality, we may assume that $y \leq z$, so $y \leq t \leq z$.
\ben
\item If $x \leq y \leq z$, then $x \leq t \leq x+2$, which is a contradiction.
\item If $y \leq x \leq z$, then $t \leq y+1 \leq x+1$ and $t \geq z-1 \geq x-1$, so $x-1 \leq t \leq x+1$, which is a contradiction.
\item If $y \leq z \leq x$, then $x-2 \leq y \leq t \leq x$, which is a contradiction.
\een
\eit
\ep

\blem
\label{lem:sc}
The triangle-square complex of $X$ is simply connected.
\elem

\bp
Assume that $\ell$ is a combinatorial loop in $X$, and fix $x \in X$ such that $x \leq \ell$. Then, for each $n \in \N$, let $\ell_n$ denote the loop $\ell \wedge (x+n)$: more precisely, if $y$ is a vertex of $\ell$, then $y \wedge (x+n)$ is a vertex of $\ell_n$. This actually defines a loop since, if $d(y,z) =1$, for instance $y \leq z \leq y+1$, then $y \wedge x \leq z \wedge x \leq (y+1) \wedge x \leq (y \wedge x)+1$, so $d(y \wedge x,z \wedge x) \leq 1$. Since also $d(y \wedge (x+n), y \wedge (x+n+1)) \leq 1$, we deduce that, for each $n \in \N$, the loops $\ell_n$ and $\ell_{n+1}$ are homotopic in the triangle-square complex of $X$.

If $N \in \N$ is such that $\ell \leq x+N$, the loop $\ell_N$ is constant equal to $\ell$, whereas $\ell_0$ is the constant loop at $x$. Hence the triangle-square complex of $X$ is simply connected.
\ep

\section{Garside categories, Garside lattices and Garside flag complexes}

\label{sec:garside}

In this section, we explicit a dictionary between categorical Garside structures and certain lattices and simplicial complexes, following Bessis (\cite{bessis2015finite}). In particular, we make connections between categorical Garside structures and the type of lattices studied in Section~\ref{sec:quotient} and Section~\ref{sec:lattice_no_quotient}.

\subsection{Definition of Garside category and an example}
Let $\mathcal C$ be a small category. One may think of $\mathcal C$ as an oriented graph, whose vertices are objects in $\mathcal C$ and oriented edges are morphisms of $\mathcal C$. Arrows in $\mathcal C$ compose like paths: $x\stackrel{f}{\to} y\stackrel{g}{\to} z$ is composed into $x\stackrel{fg}{\to} z$. For objects $x,y\in\mathcal C$, let $\mathcal C_{x\to}$ denote the collection of morphisms whose source object is $x$. Similarly we define $\mathcal C_{\to y}$ and $\mathcal C_{x\to y}$. 

For two morphisms $f$ and $g$, we define $f\preccurlyeq g$ if there exists a morphism $h$ such that $g=fh$. Define $g\succcurlyeq f$ if there exists a morphism $h$ such that $g=hf$. A nontrivial morphism $f$ which cannot be factorized into two nontrivial factors is an \emph{atom}.

The category $\mathcal C$ is \emph{cancellative} if, whenever a relation $afb=agb$ holds between composed morphisms, it implies $f=g$. $\mathcal C$ is \emph{homogeneous} if there exists a length function $l$ from the set
of $\mathcal C$-morphisms to $\mathbb Z_{\ge 0}$ such that $l(fg) = l(f) + l(g)$ and $(l(f) = 0)\Leftrightarrow$ ($f$ is a unit). If $\mathcal C$ is homogeneous, then Then $(\mathcal C_{x\to},\preccurlyeq)$ and $(\mathcal C_{\to y},\succcurlyeq)$ are posets.

We consider the triple $(\mathcal C,\mathcal C\stackrel{\phi}{\to}\mathcal C,1_\mathcal C \stackrel{\Delta}{\Rightarrow}\phi)$ where $\phi$ is an automorphism of $\mathcal C$ and $\Delta$ is a natural transformation from the identity function to $\phi$. For an object $x\in \mathcal C$, $\Delta$ gives morphisms $x\stackrel{\Delta(x)}{\longrightarrow} \phi(x)$ and $\phi^{-1}(x)\stackrel{\Delta(\phi^{-1}(x))}{\longrightarrow} x$. We denote the first morphism by $\Delta_x$ and the second morphism by $\Delta^x$. A morphism $x\stackrel{f}{\to} y$ is \emph{simple} if there exists a morphism $y\stackrel{f^\ast}{\to} \phi(x)$ such that $f f^\ast=\Delta_x$. When $\mathcal C$ is cancellative, such $f^\ast$ is unique.

\begin{defi}[\cite{bessis2006garside}]
	\label{def:Garside}
	A \emph{homogeneous categorical Garside structure} is a triple $(\mathcal C,\mathcal C\stackrel{\phi}{\to}\mathcal C,1_\mathcal C \stackrel{\Delta}{\Rightarrow}\phi)$ such that:
	\begin{enumerate}
		\item $\phi$ is an automorphism of $\mathcal C$ and $\Delta$ is a natural transformation from the identity function to $\phi$;
		\item $\mathcal C$ is homogeneous and cancellative;
		\item all atoms of $\mathcal C$ are simple;
		\item for any object $x$, $\mathcal C_{x\to}$ and $\mathcal C_{\to x}$ are lattices.
	\end{enumerate}
	It has \emph{finite type} if the collection of simple morphisms of $\mathcal C$ is finite. 
\end{defi}

A fundamental property of $\mathcal C$ is that the natural map $\mathcal C\to\mathcal G$ is an embedding, where $\mathcal G$ denotes the enveloping groupoid, as follows from the discussion in \cite[Section 2]{bessis2006garside}. For the convenience of the reader, we will remind below the existence of a normal form (see also~\cite{garside} for more details on normal forms in Garside theory, and also~\cite{bestvina_artin}).

\bpro
Let $(\mathcal C,\mathcal C\stackrel{\phi}{\to}\mathcal C,1_\mathcal C \stackrel{\Delta}{\Rightarrow}\phi)$ denote a homogeneous categorical Garside structure, and let $\mathcal G$ denote the enveloping groupoid. Then any $f \in \mathcal G$ has a unique \emph{(left greedy) normal form}, i.e. a unique way to write $f$ as a product
$$f=s_1s_2 \dots s_\ell \Delta^k,$$
where $s_1,s_2,\dots,s_\ell$ are simple elements of $\mathcal C$ with sources $x_1,x_2,\dots,x_\ell$, $k \in \Z$, $s_1 < \Delta_{x_1}$, and for all $1 \leq i \leq \ell-1$ we have the \emph{left-weighted condition}
$$s_i = s_is_{i+1} \wedge \Delta_{x_i}.$$
\epro

\bp
Let us first prove that there exists such a product, without restriction on the $s_i$'s. Since simple elements generate $\mathcal G$, we deduce that $f$ may be written as $f=s_1^{\pm 1}s_2^{\pm 1} \dots s_\ell^{\pm 1}$, where for each $1 \leq i \leq \ell$ the element $s_i$ is simple. For each simple element $s$ with source $x$ and terminal object $y$, by definition, there exists a simple element $s^\ast$ such that $s s^\ast=\Delta_x$. As a consequence, we may rewrite $s^{-1}={s^\ast}^{-1} {\Delta_x}^{-1}$. Furthermore, given any simple element $s$ with source $x$, we have $\Delta^x s = \Delta^{-1}(s) \Delta^y$, with $\Delta^{-1}(s)$ simple. Hence we see from this two properties that one can write $f$ as a product $f=s_1s_2 \dots s_\ell \Delta^k$, where $s_1,\dots,s_\ell$ are simple. Moreover, we may assume that $s_1<\Delta_{x_1}$.

\mk

We will now prove that we can assume the left-weighted condition: assume that $1 \leq i \leq \ell-1$ is such that $s_i < s'_i=s_is_{i+1} \wedge \Delta_{x_i}$. Then $s'_i$ is simple, and we may write $s_is_{i+1}=s'_is'_{i+1}$, we will show that $s'_{i+1}$ is simple. Let us also write $s'_i=s_it$: since $\mathcal C$ is cancellative, we have $ts'_{i+1}=s_{i+1}$. Since $s_{i+1}$ is simple, we have $s_{i+1} s_{i+1}^\ast = \Delta_{x_{i+1}}$, so $ts'_{i+1}s_{i+1}^\ast = \Delta_{x_{i+1}}$. In particular, $t^\ast=s'_{i+1}s_{i+1}^\ast$ is simple, so $s'_{i+1}$ is simple. Hence there exists a left-weighted normal form for $f$.

\mk

Let us now prove that this form is unique. Without loss of generality, we may assume that $f \in \mathcal C$. Let us consider a left-weighted forms $f=s_1s_2 \dots s_\ell \Delta^k$, with $k \geq 0$. Then it is easy to see by induction that, for all $1 \leq i \leq \ell$, we have $s_1s_2 \dots s_i = f \wedge \Delta^i_{x_1}$. Uniqueness follows.
\ep

\begin{defi}
	A \emph{Garside category} is a category $\mathcal C$ that can be equipped with $\phi$ and $\Delta$ to obtain a homogeneous categorical Garside structure. A \emph{Garside groupoid} is the enveloping groupoid of a Garside category. Informally speaking, it is a groupoid obtained by adding formal inverses to all morphisms in a Garside category.
	
	Let $x$ be an object in a groupoid $\mathcal G$. The \emph{isotropy} group $\mathcal G_x$ at $x$ is the group of morphisms from $x$ to itself. A \emph{weak Garside group} is a group isomorphic to the isotropy group of an object in a Garside groupoid. 
	
	A \emph{Garside monoid} is a Garside category with a single object and a \emph{Garside group} is a Garside groupoid with a single object.
\end{defi}

\begin{exe}
	Let $\mathcal A$ be a finite central arrangement in $\R^n$, i.e.\ a finite collection of linear hyperplanes in $\R^n$. Let $ch(\mathcal A)$ be the set of chambers (connected components of the complement of the hyperplanes in $\mathcal A$). We consider an oriented graph $\Gamma$, whose vertices are in 1-1 correspondence with the collection of chambers, and we draw a pair of oriented edges going in opposite directions between two vertices if the associated chambers are adjacent along a hyperplane.

	A \emph{positive path} on $\Gamma$ is an edge path from one vertex to another vertex which goes along the positive orientation on each edge. Take a positive path $f_1$ and a subpath $g$ of $f_1$ which is a geodesic between its endpoints with respect to the path metric on $\Gamma$. An \emph{elementary homotopy} of $f_1$ is the procedure of replacing the subpath $g$ of $f_1$ by another positive subpath which is a geodesic between the two endpoints of $g$. Two positive paths are \emph{equivalent} if they differ by a finite sequence of elementary homotopies.
	
	Now we consider the category $\mathcal C$ as follows. Objects of $\mathcal C$ are vertices of $\Gamma$, and morphisms of $\mathcal C$ are equivalence classes of positive paths from one vertex to another vertex. There is an orientation-preserving automorphism $\alpha:\Gamma\to\Gamma$ arising from the central symmetry of $\R^n$ with respect to the origin. Note that $\alpha$ induces an automorphism of the category $\mathcal C\stackrel{\phi}{\to}\mathcal C$, which is a functor sending object $x$ to $\alpha(x)$, and the morphism represented by a path $f$ to the morphism represented by $\alpha(f)$.
	For each object $x\in \mathcal C$, let $\Delta_x$ be the morphism from $x$ to $\phi(x)$ represented by a positive geodesic in $\Gamma$ from $x$ to $\phi(x)$. One readily verifies that the family of morphisms $\{\Delta_x\}_{x\in\mathrm{Obj}(\mathcal C)}$ gives a natural transformation between the identity function and the functor $\phi$, i.e. for each morphism $[f]$ in $\mathcal C$ represented by a positive path $f$ from $x$ to $y$, we have $\Delta_x\phi([f])=[f]\Delta_y$.

	It is shown in \cite{deligne1972immeubles} (see also \cite[Example 3.4]{bessis2006garside}) if the arrangement $\cal A$ is \emph{simplicial}, namely, hyperplanes in $\cal A$ cuts the unit sphere of $\R^n$ into a simplicial complex, then $\mathcal C$ is a homogeneous categorical Garside structure in the above sense (the non-trivial part is property (4) of Definition~\ref{def:Garside}). Moreover, the associated weak Garside group is isomorphic to the fundamental group of the complement of the complexification of hyperplanes of $\mathcal A$ in $\mathbb C^n$.
\end{exe}

\subsection{From Garside category to Garside lattice}

The type of lattices studied in Section~\ref{sec:quotient} and Section~\ref{sec:lattice_no_quotient}, if they are assumed homogeneous, are what we call a Garside lattice, which we define now.

\bdf
A \emph{Garside lattice} is a pair $(L,\varphi)$, where $L$ is a homogeneous lattice and $\varphi$ is an increasing automorphism of $L$, such that, for any $x,y \in L$, there exists $k \in \N$ such that $x \leq \varphi^k(y)$.
\edf

We now see that a categorical Garside structure naturally gives rise to a Garside lattice.

\bpro \label{pro:Garside_lattice_from_category}
Let $(\mathcal C, \phi,\Delta)$ be a categorical Garside structure, and let $\mathcal G$ be the associated Garside groupoid. Let $x$ be an object in $\mathcal G$, and let $L_x$ denote the set of morphisms from $x$ in $\mathcal G$. Let us consider the map $\psi: f \in L_x \mapsto \Delta_x \phi(f)$ of $L_x$. We endow $L_x$ with the partial order $\le$ such that $g\le h$ if $h=gf$ with $f\in\mathcal C$. Then $(L_x,\le)$ is a lattice such that 
\begin{enumerate}
	\item any non-empty upper bounded set in $L_x$ has a join;
	\item $\psi$ is an increasing automorphism of $(L_x,\le)$;
	\item for any $g,h \in L_x$, there exists $k \in \N$ such that $\psi^{-k} (g) \leq h \leq \psi^{k} (g)$.
\end{enumerate}
\epro

\bp
Recall that the natural map $\mathcal C\to \mathcal G$ is an embedding. Given $g\in L_x$, by Definition~\ref{def:Garside} (4), the collection $(L_x)_{\ge g}$ of all elements of $L_x$ that is $\ge g$ form a lattice under the order $\le$. To prove that $(L_x,\le)$ is a lattice, we will show that any two elements in $L_x$ have a lower bound. Indeed, given $f,g\in L_x$, by Definition~\ref{def:Garside} (2) and (3), we can write $f=gs^{\eps_1}_1s^{\eps_2}_2\cdots s^{\eps_n}_n$, where $\eps_i=\pm 1$, and each $s_i\in \mathcal C$ is a simple element. We claim $\psi(f)$ is of the form $gh(s')^{\eps_2}_2\cdots (s')^{\eps_n}_n$ where $h\in\mathcal C$ and each $s'_i\in\mathcal C$ is simple. Indeed, as $\Delta$ gives a natural transformation from the identity function to $\phi$, we know $\psi(f)=\Delta_x\phi(g)\phi(s^{\eps_1}_1)\cdots \phi(s^{\eps_n}_n)=g\Delta_{t(g)}\phi(s^{\eps_1}_1)\cdots \phi(s^{\eps_n}_n)$ where $t(g)$ denotes the terminal object of $g$. Then the claim is clear if $\eps_1=1$. If $\eps_1=-1$, we find $s^*_1$ such that $s_1s^*_1=\Delta_{t(s^{\eps_1}_1)}$. As $\Delta$ is a natural transformation, we know $\Delta_{t(g)}\phi(s^{\eps_1}_1)=s^{\eps_1}_1\Delta_{t(s^{\eps_1}_1)}$. Let $y$ be the starting point of $s_1$. Then $\Delta_{t(s^{\eps_1}_1)}=\Delta_y=s_1s^*_1$ for some $s^*_1\in\mathcal C$. Thus $\Delta_{t(g)}\phi(s^{\eps_1}_1)=s^*_1$ and the claim follows. By applying the claim several times, we know $\psi^n(f)=gh$ with $h\in\mathcal C$. Thus $g\le \psi^n(f)$. Similarly $f\le \psi^n(g)$. Thus property (3) of the proposition holds. In particular, any two elements in $L_x$ have a lower bound. Thus $L_x$ is a lattice. As $L_x$ is homogeneous, chains in intervals have bounded length. Hence, given any non-empty subset $A \subset L_x$ with an upper bound $f \in L_x$, by considering the joins of increasing finite subsets of $A$, one sees that $A$ itself has a join. Hence property (1) holds.

Note that, if $f \in \mathcal G$ has source $x$ and terminal object $y$, then $\psi(f)=\Delta_x \phi(f) = f \Delta_y$. From this, one sees that $\psi$ is invertible and increasing, hence property (2) holds.
\ep

\subsection{Garside categories, Garside lattices and Garside flag complexes}

\label{subsec:dictionary}

We will relate three Garside notions: categorical Garside structures, Garside lattices, and Garside flag complexes, which we define now. These are flag simplicial complexes with a local lattice structure and a special automorphism.

\bdf
\label{def:garside flag}
A \emph{Garside flag complex} is a pair $(X,\varphi)$, where $X$ is a simply connected flag simplicial complex with a consistent total order on each simplex, a labeling of edges $\ell:EX\to \mathbb Z_{>0}$ and $\varphi$ is an order-preserving and $\ell$-preserving automorphism of $X$. Let $<$ be the binary relation on the set of vertices of $X$ such that $x<y$ if $x$ and $y$ are adjacent in $X$ and $x<y$ with respect to the order on the edge connecting $x$ and $y$. We also require  the following:
\ben
\item for any 2-dimensional simplex with vertices $a<b<c$, we have $\ell(ab)+\ell(bc)=\ell(ac)$;
\item $a\le b$ if and only if $b\le \varphi(a)$;
%\item For any simplex $\sigma$ of $X$, we have that $\sigma \cup \varphi(\min \sigma)$ is a simplex of $X$.
\item for any $x\in X$, the set $[x,\varphi(x)]=\{z\in X^{(0)}\mid x\le z\ \textrm{and}\ z\le \varphi(x)\}$ with the binary relation $<$ is a lattice.
\een
\edf
%
%We say that a group $G$ acts by automorphisms on a categorical Garside structure $(\mathcal C, \phi,\Delta)$ if $G$ acts by automorphisms on the category $\mathcal C$ and its action is compatible with $\phi$ and $\Delta$ in the following sense:
%\begin{enumerate}
%	\item the action of $G$ commutes with the action of $\phi$;
%	\item $\Delta$ is a natural transformation between the identity functor and each functor in $G$.
%\end{enumerate}
%If a group $G$ acts by automorphisms on a categorical Garside structure $(\mathcal C, \phi,\Delta)$ such that the restriction on the action to the set of objects is free, then one may consider the quotient categorical 
%Garside structure, whose set of objects is the quotient of the set of objects of $\mathcal C$ by $G$. These two categorical Garside structures, however, have the same set of morphisms starting a given object. We will call two categorical Garside structures \emph{equivalent} if they are quotients of the same categorical Garside structure by actions of automorphisms which are free on the set of objects. 

%If $(\mathcal C,\phi,\Delta)$ is a connected categorical Garside structure, then $(\mathcal C',\phi',\Delta')$ is connected categorical Garside structure where $\phi'$

%\begin{lem}
%If $(\mathcal C,\phi,\Delta)$ is a connected categorical Garside structure, then $(\mathcal C',\phi',\Delta')$ is connected categorical Garside structure where $\phi'$
%\end{lem}

\mk

Note that, given any Garside flag complex $(X,\varphi)$ and for any simplex $\sigma$ in $X$, we have that $\sigma \cup \varphi(\min \sigma)$ is a simplex, whose maximal element is $\varphi(\min \sigma)$.

\mk

Theses notions enable us to write the following dictionary between categorical Garside structures, Garside lattices and Garside flag complexes. 

%A category is \emph{special}, if for any pair of objects $x,y$, there is at most one morphism from $x$ to $y$.
A categorical Garside structure $(\mathcal C, \phi,\Delta)$ is \emph{special} if 
\ben
\item for any pair of objects $x,y$, there is at most one morphism from $x$ to $y$;
\item the nerve of $\mathcal C$ (cf. \cite[Section 7]{bessis2006garside}) is simply-connected.
\een

Each connected categorical Garside structure $(\mathcal C, \phi,\Delta)$ gives a canonical special categorical Garside structure via a ``universal covering'' construction as in \cite[Section 6]{bessis2006garside} as follows. Let $\mathcal G$ be the associated Garside groupoid. 
Take an object $x$ of $\mathcal G$. 
Consider the category $\mathcal C'$ whose objects corresponds to morphisms of $\mathcal G$ starting at $x$. 
Given two objects $x'_f$ and $x'_g$ of $\mathcal C'$, corresponding to two morphisms $f,g$ of $\mathcal G$ starting at $x$, we assign a morphism from $x'_f$ to $x'_g$ if $g=fh$ for $h\in \mathcal C$. One readily verifies that $\mathcal C'$ is indeed special. Let $\phi'$ be the automorphism of $\mathcal C'$ sending the morphism $f$ of $\mathcal G$ from $x$ to $y$ to the morphism $f\Delta_y$ of $\mathcal G$. Now we check that $\phi'$ is indeed an automorphism of the category $\mathcal C'$. This amounts to check that, if the elements $f,g$ of $\mathcal G$ start at $x$ and end at $y$ and $z$ respectively, and $f=gh$ for $h\in \mathcal C$, then $f\Delta_y=g\Delta_zh'$ for $h'\in \mathcal C$. Indeed, $f\Delta_y=gh\Delta_y=g\Delta_zh'$ as $\Delta$ is a natural transformation from the identity functor on $\mathcal C$ to $\phi$. For each object $x'_f$ of $\mathcal C'$ associated with the morphism $f$ of $\mathcal C$ from $x$ to $y$, let $\Delta'_{x'_f}$ be the unique morphism of $\mathcal C'$ from $x'_f$ to $x'_{f\Delta_y}$. Then one readily verifies that $(\mathcal C', \phi',\Delta')$ gives a special categorical Garside structure (see also \cite[Lemma 6.2]{bessis2006garside}). As the category $\mathcal C$ is connected, $(\mathcal C', \phi',\Delta')$ does not depend on the choice of the object $x$ of $\mathcal G$, as different choices of $x$ give isomorphic special categorical Garside structure.

\bthm \label{thm:dictionary_garside}
\label{theorem:dictionary}
The following objects are in one to one correspondence:
\ben
\item A connected special categorical Garside structure, up to isomorphism.
\item A Garside lattice, up to isomorphism.
\item A Garside flag complex, up to isomorphism.
\een
\ethm

\bp
First we establish the correspondence between item 1 and item 2.
Given a connected special categorical Garside structure $(\mathcal C, \phi,\Delta)$, let $\mathcal G$ be the associated Garside groupoid. Note that $\mathcal G$ and $\mathcal C$ have the same set of objects.
Fix an object $x $ of $ \mathcal G$ and let $L$ be the collection of morphisms of $\mathcal G$ starting at $x$. Let $\varphi:L\to L$ sending a morphism $f$ of $\mathcal G$ from $x$ to $y$ to the morphism $f\Delta_y$. We can define an order on $L$ by $f\le g$ for morphisms $f,g$ of $\mathcal G$ starting from $x$ if and only if $f=gh$ for $h$ being a morphism of $\mathcal C$. Note that the homogeneous structure on $\mathcal C$ implies that $(L,\le)$ is indeed a poset, and it is homogeneous. By the discussion of normal forms in \cite[Section 2]{bessis2006garside} that every two elements in $(L,\le)$ have an upper bound and a lower bound, now item 4 of Definition~\ref{def:Garside} implies that $(L,\le)$ is a lattice. It follows that $(L,\varphi)$ is a Garside lattice. As $\mathcal C$ is connected, the isomorphism type of this Garside lattice does not depend on the choice of $x$.

%The length function on the edges of $\operatorname{gar}(\mathcal G,S,x)$ comes from the length function on morphisms of $\mathcal C$. Note that, since $\mathcal G$ is connected, $\operatorname{gar}(\mathcal G,S,x)$ does not depend on $x$. 

Given a Garside lattice $(L,\varphi)$, one may consider the category $\mathcal C'$ whose objects are elements of $L$, and there is a unique morphism from $x$ to $y$ if and only if $x \leq y$. Then $\mathcal C'$ is special, where the simply-connectedness requirement follows from the lattice property (for a given loop in the nerve of $\mathcal C'$, we can find a common lower bound $x_0\in L$ for the vertex set of this loop, and use the lattice property to deform this loop towards $x_0$ in the nerve). The automorphism $\varphi$ of $L$ induces an automorphism $\phi'$ of $\mathcal C'$. Given any morphism $x<y$, one may also define $\phi'(x<y)$ as the unique morphism $\varphi(x)<\varphi(y)$.
 Given any $x \in L$, one may define $\Delta'(x)$ as the unique morphism $x<\varphi(x)$.  So $(\mathcal C', \phi',\Delta')$ is a special categorical Garside structure. 
 
 Now we check the procedures in the previous paragraphs are inverses of each other. We will only verify one direction, and leave the other direction to the reader.
Consider $(\mathcal C,\phi,\Delta)\to (L,\varphi)\to (\mathcal C',\phi',\Delta')$ and we show the two Garside categories $\mathcal C$ and $\mathcal C'$ are isomorphic. By construction, objects of $\mathcal C'$ are in 1-1 correspondence with the set of morphisms $\mathcal G_{x\to}$ of $\mathcal G$ starting at $x$. As $\mathcal C$ is connected, for any object $y$ of $\mathcal G$, there is a morphism in $\mathcal G_{x\to}$ ending at $y$.
As any morphism of $\mathcal G$ from $x$ to $y$ is a finite composition of morphisms in $\mathcal C$ or their inverses, the specialness of $\mathcal C$ implies that there is a unique morphism in $\mathcal G$ from $x$ to $y$. Thus objects of $\mathcal C'$ are in 1-1 correspondence with objects in $\mathcal C$ by sending an element in $\mathcal G_{x\to}$ to its end object. Now one readily checks such 1-1 correspondence actually induces isomorphism an of the associated categorical Garside structure. 
 
%Now take a special categorical Garside structure $(\mathcal C,\phi,\Delta)$ and produce the Garside lattice $(L,\varphi)$ with $L$ being the collection of morphisms of $\mathcal G$ starting at an object $x$ of $\mathcal G$.
 
%One readily checks that the procedure in this paragraph and the procedure in the previous paragraph are inverses of each other.

 %, which gives a Garside groupoid $\mathcal G'$. Let $x'$ denote the object of $\mathcal G'$ consisting of the identity morphism $x \leq x$. Let $G'$ denote the weak Garside group that is the isotropy group of $\mathcal G'$ at $x'$.

%We also see that the lattice from Proposition~\ref{pro:Garside_lattice_from_category} consisting of morphisms from $x$ forms a Garside lattice.

\mk

Now we establish the correspondence between item 2 and item 3.
Given a Garside flag complex $(X,\varphi)$, we consider the binary relation $<$ on the vertex set $L$ of $X$. Suppose $a,b,c,d\in L$ such that $a\le b, b\le c$, $c\le d$. Now we verify condition $(\ast)$ in Definition~\ref{defi:weak_order}. First assume $(b,c,d)$ and $(a,b,d)$ are transitive triples. Then $a\le d$, $b\le d$ and $c\le d$. Then Definition~\ref{def:garside flag} (2) implies that $\{a,b,c,d\}\subset [\varphi^{-1}(d),d]$. As $[\varphi^{-1}(d),d]$ is a lattice by Definition~\ref{def:garside flag} (3), in particular it is a poset, hence $(a,b,c)$ and $(a,c,d)$ are transitive triples. Conversely, if $(a,b,c)$ and $(a,c,d)$ are transitive triples, then $a\le b, a\le c, a\le d$, implying $\{a,b,c,d\}\subset [a,\varphi(a)]$ by Definition~\ref{def:garside flag} (2). Now we deduce from Definition~\ref{def:garside flag} (3) as before that $(b,c,d)$ and $(a,b,d)$ are transitive triples. Thus the binary relation gives a weak order on $L$. This weak order is homogeneous by Definition~\ref{def:garside flag} (1). The map $\varphi:L\to L$ is an automorphism of this weak order, which satisfies the assumptions of Theorem~\ref{thm:lattice ltg} by the definition of a Garside flag complex. Thus for any $x \in L$, the posets $L_{\leq x}$ and $L_{\geq x}$ are lattices. Since $X$ is connected, we deduce that, for any $x,y \in L$, there exists $k \in \N$ such that $x \leq \varphi^k(y)$. In particular, one sees that any two elements of $L$ have a lower bound, so $L$ is a lattice. Hence $L$ generates a Garside lattice as in Theorem~\ref{thm:lattice ltg}.

Conversely, consider a Garside lattice $(L,\varphi)$. Let $X$ denote the simplicial complex whose vertices are elements in $L$, with a $k$-simplex for each length $k$ chain $x_0<x_1<\dots<x_k$ in $L$ such that $x_k \leq \varphi(x_0)$. The automorphism $\varphi$ of $L$ induces an order-preserving automorphism (still denoted $\varphi$) of $X$. Since $L$ is homogeneous, there exists a function $\ell$ from the set of comparable pairs of $L$ to $\Z_{>0}$. This function restricts to $\ell:E(X) \ra \Z_{>0}$ such that, for any $2$-dimensional simplex with vertices $a<b<c$, we have $\ell(ab)+\ell(bc)=\ell(ac)$. Moreover, if $a<b$ are the vertices of an edge of $X$, then by definition we have $b \leq \varphi(a)$. Conversely, if $b<\varphi(a)$ are the vertices of an edge of $X$, then $\varphi(a) \leq \varphi(b)$, hence $a \leq b$. Finally, for any vertex $x$ in $X$, the set $[x,\varphi(x)]$ is an interval in the lattice $L$, hence it is a lattice.
Moreover, this procedure of passing from Garside lattice to Garside flat complex and the procedure in the previous paragraph are inverses of each other, which gives the correspondence between item 2 and item 3.
%
%\mk
%
%We will check that the composition of these three constructions, starting from a categorical Garside structure $(\mathcal C, \phi,\Delta)$, gives an equivalent categorical Garside structure $(\mathcal C', \phi',\Delta')$. Note that the weak Garside group $G$ acts on the objects of $\mathcal C'$ by precomposition, so $G$ acts by automorphisms on $\mathcal C'$. Since $\phi'$ and $\Delta'$ are defined by postcomposition, we deduce that the action of $G$ is compatible with $\phi'$ and $\Delta'$. Hence $G$ acts by automorphisms on $(\mathcal C', \phi',\Delta')$.
%
%The set of objects of the quotient $\mathcal C' /G$ then identifies canonically with the set of objects of $\mathcal C$, and the set of morphisms of $\mathcal C' /G$ identifies canonically with the set of morphisms of $\mathcal C$. The quotient of $\phi'$ and $\Delta'$ by $G$ also identify canonically with $\phi$ and $\Delta$. Hence the categorical Garside structures $(\mathcal C, \phi,\Delta)$ and $(\mathcal C', \phi',\Delta')$ are equivalent.
\ep

\begin{rmk}
Here is a way to go directly from item 1 of the above theorem to item 3.	Given a special categorical Garside structure $(\mathcal C, \phi,\Delta)$, let $\mathcal G$ denote the associated Garside groupoid, fix an object $x \in \mathcal G$, and let $G$ denote the weak Garside group that is the isotropy group of $\mathcal G$ at $x$. Then the simplicial complex $\operatorname{gar}(\mathcal G,S,x)$ from~\cite[Definition~B.11]{bessis2015finite}, whose vertex set is the set of morphisms from $x$ and whose edges correspond to morphisms of $\mathcal G$ from $x$ that differ by a simple morphism in $\mathcal C$, is a Garside flag complex. More precisely, the map $\varphi$ in Definition~\ref{def:garside flag} sends a morphism $f$ of $\mathcal G$ from $x$ to $y$ to the morphism $f\Delta_y$. The length function on the edges of $\operatorname{gar}(\mathcal G,S,x)$ comes from the length function on morphisms of $\mathcal C$.
	Note that, since $\mathcal G$ is connected, $\operatorname{gar}(\mathcal G,S,x)$ does not depend on $x$. 
\end{rmk}

We also have a similar characterization of (weak) Garside groups.

\bthm \label{thm:dictionary_garside_group}
A group $G$ is a Garside group (resp. weak Garside group) if and only if there exists a Garside flag complex $(X,\varphi)$ (or equivalently on a Garside lattice) such that $G$ can be realized as a group of order-preserving automorphisms of $X$ commuting with $\varphi$, acting freely and transitively (resp. freely) on vertices of $X$.

Moreover, the group $G$ is a (weak) Garside group of finite type if and only if $X$ can be chosen such that the action of $G$ is cocompact.
\ethm

\bp
If $G$ is a (weak) Garside group, it is clear from Theorem~\ref{thm:dictionary_garside} that $G$ acts on the associated Garside lattice satisfying these properties.

\mk

Conversely, assume that $G$ acts on a Garside lattice $X$ satisfying these properties. Then, if we fix a vertex $x \in X$, we may consider the categorical Garside structure whose objects are elements of $X$ bigger than $x$ as in Theorem~\ref{thm:dictionary_garside}. Its quotient by $G$ is a categorical Garside structure $(\mathcal C,\phi, \Delta)$, and the isotropy group at the image of $x$ in the associated groupoid coincides with $G$.
\ep

\bexe
As a very simple example, consider the $n$-strand braid group $B_n$, with the standard Garside element $\Delta_n$ and positive monoid $B_n^+$. Note that $B_n^+$ plays the role of the categorical Garside structure (with only one object), and $B_n$ is the associated Garside groupoid. 
Let $L$ denote the lattice consisting of elements in $B_n$, where $g \leq h$ if and only if $h \in gB_n^+$. The automorphism $\varphi: g \mapsto g\Delta$ of $L$ turns $L$ into a Garside lattice. Now the group $B_n$ acts freely transitively on $L$ by left multiplication, commuting with $\varphi$, so we recover that $B_n$ is a Garside group. Moreover, any subgroup of $B_n$ is a weak Garside group, and any finite index subgroup of $B_n$ is a weak Garside group of finite type. This is the case of the pure braid group, in which case the categorical Garside structure coincides with the structure described in the example about simplicial hyperplane arrangements.
\eexe
\section{Applications}
\subsection{Euclidean buildings}

\bthm \label{thm:euclidean_building_wm}
Let $X$ denote a Euclidean building of type $\tilde{A}_n$. Then the $1$-skeleton of $X$ is a weakly modular graph.
\ethm

\bp
Let $m=n+1$. Let us denote the type function $\tau : X^{(0)} \ra \Z/m\Z$ (see for instance~\cite[Section~6.9]{abramenko_brown}). Let us consider
$$L=\{(x,k) \in X^{(0)} \times \Z \st \tau(x) \equiv k [m]\}.$$
Let us consider the order relation on $L$ generated by $(x,k) \leq (x',k')$ if $x$ and $x'$ are adjacent in $X$ and $k \leq k'$. According to~\cite{hirai_uniform_modular}, $L$ is a lattice.

\mk

Note that this lattice is quite explicit in the Bruhat-Tits case: let $X$ denote the Bruhat-Tits building of $\PGL(n,\K)$, where $\K$ is a non-Archimedean local field. The vertex set of $X$ may be described as the set of homothety classes of ultrametric norms on $\K^n$. The extended Bruhat-Tits building $\hat{X}$ of $\GL(n,\K)$ may be described as the set of ultrametric norms on $\K^n$. There is a canonical simplicial map $\hat{X} \ra X$, such that the preimage of a vertex is isomorphic to $\Z$. The vertex set $L$ of $\hat{X}$, with the natural order on norms on $\K^n$, is in fact a lattice (see also~\cite{haettel_injective_buildings})

\mk

Let us consider the action of $\Z$ on $L$ by $p \cdot (x,k)=(x,k+pm)$. This action satisfies the assumptions of Theorem~\ref{thm:quotient wm}, so the quotient graph of $L/\Z$, which coincides with $X^{(1)}$, is weakly modular.
\ep

In fact, we can see Euclidean buildings of type $\tilde{A}_n$ as categorical Garside structures, as follows.

\bthm
\label{thm:building}
Let $X$ denote a Euclidean building of type $\tilde{A}_n$, and let $G$ denote an automorphism group of $X$. Then $X^{(0)}/G$ is the set of objects of a categorical Garside structure. Moreover, for any subgroup $H$ of $G$ acting freely on $X$, the group $H \times \Z$ is weakly Garside. 
\ethm

\bp
According to the proof of Theorem~\ref{thm:euclidean_building_wm}, the lattice $L \subset X^{(0)} \times \Z/n\Z$ is a Garside lattice. So according to Theorem~\ref{thm:dictionary_garside}, we deduce that $X^{(0)}$ is the set of objects of a categorical Garside structure on which $G$ acts by automorphisms. Hence $X^{(0)}/G$ is the set of objects of a categorical Garside structure.

\mk

Moreover, given any subgroup $H$ of $G$ acting freely on $X$, the group $H \times \Z$ acts freely by automorphisms on $L$, so according to Theorem~\ref{thm:dictionary_garside_group} it is a weak Garside group.
\ep

\bexe
For instance, the Bruhat-Tits building of $\PGL(n,\K)$, for a non-Archimedean valued field $\K$, may be described as the categorical Garside structure associated to the following Garside germ $(\mathcal C,S)$. The category $\mathcal C$ has one object, and the set of simple morphisms $S$ coincide with the poset of vector subspaces of $\k^n$, where $\k$ is the residue field. The partial composition of morphisms coincides with the natural inclusion order on $S$. If one considers the Bruhat-Tits building of $\SL(n,\K)$ instead, one may naturally consider a similar Garside structure, whose objects are $\Z/n\Z$.
\eexe

\brk
In the view of~\cite{haettel_injective_buildings}, one may also consider the space $X$ of all convex symmetric bodies of $\R^n$, which can be understood as the injective hull of the symmetric space $\GL(n,\R)/O(n)$, as a "continuous" categorical structure, where the homogeneous assumption should be replaced with a graded with values in $\R$. Hence torsion-free subgroups of $\GL(n,\R)$ can therefore be described as "continuous" weak Garside groups.
\erk

\subsection{Exotic Garside groups}

In fact, we can use the previous result to produce actual Garside groups with exotic properties.

\mk

We will describe the construction of Theorem~\ref{thm:euclidean_building_wm} in the more general and simple framework of a \emph{typed $\tilde{A}_2$ complex}: it is a connected, simply connected CAT(0) equilateral triangle complex $X$, with a type function $\tau:X^{(0)} \ra \Z/3\Z$ such that adjacent vertices have different types. Let us say that an automorphism $g$ of $X$ is type-rotating if there exists $a \in \Z/3\Z$ such that $\tau \circ g = \tau + a$.

\bthm \label{thm:product_Z_Garside}
Let $X$ denote a locally finite typed $\tilde{A}_2$ complex, and let $G$ denote a group of type-rotating automorphisms of $X$ acting freely and transitively on vertices of $X$. Then a (virtually split) central extension of $G$ by $\Z$ is a Garside group.
\ethm

\bp
Let $Y$ denote the flag simplicial complex with vertex set $\{(x,n) \in X^{(0)} \times \Z \st \tau(x) \equiv n [3]\}$, with an ordered edge from $(x,n)$ to $(y,m)$ if and only if:
\bit
\item either $x$ is adjacent to $y$, and $n<m<n+3$,
\item or $x=y$, and $m=n+3$.
\eit
Note that, in case $X$ is an $\tilde{A}_2$ building, the complex $Y$ is called a building of extended type $\tilde{A}_2$. Consider the map $\ell=EY \ra \Z_{>0}$ sending the edge between $(x,n)$ and $(y,m)$ to $m-n$.  Consider the automorphism $\varphi:(x,n) \mapsto (x,n+3)$ of $Y$. We see easily that $Y$ is a Garside flag complex.

\mk

Let $H$ denote the group of automorphisms of $Y$ commuting with $\varphi$, whose induced action on the quotient $X$ coincides with $G$. Note that since we assume the action of $G$ is type rotating, the natural homomorphism $H\to G$ is surjective. Then the group $H$ is a central extension of $G$ by $\Z$. It is easy to see that $H$ acts freely and transitively on vertices of $Y$, and commuting with the automorphism $\varphi$. According to Theorem~\ref{thm:dictionary_garside_group}, we see that $H$ is a Garside group.
\ep

\bcor \label{cor:exotic_garside}
There exists a non-linear finite type Garside group $G_1$.

There exists a finite type Garside group $G_2$ without any non-trivial action on a Gromov-hyperbolic space. Moreover, any normal subgroup of $G_2$ is either finite, or has virtually cyclic quotient.

In addition, for each of these Garside groups $G_i$, the central quotient $G_i/Z(G_i)$ has Lafforgue's strong Property (T).
\ecor

\bp
According to~\cite[Theorem~1.1]{radu_exotic_lattice}, there exists an exotic locally finite $\tilde{A}_2$ building $X_1$ with a group $H_1$ of type-rotating automorphisms of $X_1$ acting freely and transitively on vertices of $X_1$. According to Theorem~\ref{thm:product_Z_Garside}, there exists a central extension $G_1$ of $H_1$ by $\Z$ which is a finite type Garside group. According to~\cite{bader_caprace_lecureux}, the group $H_1$ is not linear. Note that the group $G_1$ has an index $3$ subgroup $G'_1$ isomorphic to the split central extension $H'_1 \times \Z$, where $H'_1$ is the index $3$ type-preserving subgroup of $H_1$. Since $H'_1$ is not linear, $G'_1$ is not linear.

According to~\cite{cartwright_mantero_steger_zappa_I,cartwright_mantero_steger_zappa_II}, there are many examples of algebraic locally finite $\tilde{A}_2$ building $X_2$ with a group $H_2$ of type-rotating automorphisms of $X_2$ acting freely and transitively on vertices of $X_2$. According to Theorem~\ref{thm:product_Z_Garside}, there exists a central extension $G_2$ of $H_2$ by $\Z$ which is a finite type Garside group. According to~\cite{haettel_hyprig}, the group $H_2$ has only elliptic and parabolic action by isometries on Gromov-hyperbolic spaces. We deduce that any action of $G_2$ on a Gromov-hyperbolic space is trivial: more precisely, it is elliptic, parabolic or lineal.

Finally, according to~\cite{lecureux_delasalle_witzel_T}, the groups $H_1$ and $H_2$ satisfy Lafforgue's strong Property (T).
\ep

Note that it is believed that lattices in exotic $\tilde{A}_2$ buildings also satisfy a hyperbolic rigidity, either by adapting the proof from~\cite{haettel_hyprig}, or by using the powerful theory developed by Bader and Furman (see~\cite{bader_furman_icm}).

\subsection{Weakly modular thickenings of buildings}

Here we prove Theorem~\ref{mainthm:many_buildings_thickening_wm}, showing that Conjecture~\ref{conj:building_wm} is true for many buildings.

Let $X$ denote the $1$-skeleton of a building. We say that a graph $\Gamma$ is an equivariant thickening of $X$ if $\Gamma$ contains $X$ as a subgraph, $X$ is quasi-isometric to $\Gamma$, and the automorphism group of $X$ extends as an automorphism group of $\Gamma$. Combining our result for Euclidean buildings of type $\widetilde A_n$ with other results, we deduce the following.

\bthm \label{thm:many_buildings_thickening_wm}
The following buildings have equivariant weakly modular thickenings.
\bit
\item Any spherical building.
\item Any Euclidean building of type $\widetilde A_n$, $\widetilde B_n$, $\widetilde C_n$, $\widetilde D_n$ or $\widetilde G_2$.
\item Any right-angled building.
\item Any rank $3$ building.
\item Any Gromov-hyperbolic building.
\eit
\ethm

\bp\
\bit
\item The spherical case is trivial, since the full graph on the $X$ is a weakly modular equivariant thickening.
\item If $X$ is a Euclidean building of type $\widetilde A_n$, Theorem~\ref{thm:euclidean_building_wm} states that $X$ is weakly modular, without the need of a thickening. According to~\cite{haettel_injective_buildings}, any Euclidean building of type $\widetilde B_n$, $\widetilde C_n$ or $\widetilde D_n$ has an equivariant thickening which is Helly, thus weakly modular.
\item Any right-angled building has a $1$-skeleton which is a median graph, and is in particular weakly modular.
\item According to~\cite{przytycki_schwer_systolic_buildings}, any rank $3$ building which is not of type $(2,4,4)$, $(2,4,5)$ or $(2,5,5)$ has an equivariant thickening which is systolic, thus weakly modular. The type $(2,4,4)$ is the affine type $\widetilde{C_2}$, which is already covered. The types $(2,4,5)$ or $(2,5,5)$ are of hyperbolic type, and are covered by the last case.
\item Any Gromov-hyperbolic building has a cobounded Helly hull according to~\cite{lang}, and in particular has an equivariant thickening which is weakly modular.
\eit
\ep

\subsection{Weak Garside groups}
Let $(\mathcal C, \phi,\Delta)$ be a categorical Garside structure, let $x$ be an object in $\mathcal C$, and let $L_x$ denote the set of morphisms from $x$ in the groupoid $\mathcal G$ associated with $\mathcal C$.
\begin{defi}
	Let us consider a weak Garside group $G_x$ associated to an object $x$ in a categorical Garside structure $(\mathcal C, \phi,\Delta)$ and set of simple morphisms $S$. The \emph{weak Cayley graph} $\Cay(G_x,S)$ of $G_x$ is the graph with vertex set $L_x$, the set of morphisms from $x$, with an edge between $f$ and $g$ if there exists a simple morphism $s \in S$ such that $f=gs$ or $g=fs$. The \emph{Garside automorphism} of $\Cay(G_x,S)$ is the map $f \mapsto \Delta_x \phi(f)$. Note that the group $G_x$ acts on $\Cay(G_x,S)$ by precomposition.
\end{defi}

\bthm \label{thm:garside_cayley_wm}
Let $G_x$ denote a weak Garside group, with set of simple morphisms $S$. Then the weak Cayley graph $\Cay(G_x,S)$ and its quotient $\Cay(G_x,S) / \<\Delta_x\>$ are weakly modular graphs.
\ethm

\bp
By Proposition~\ref{pro:Garside_lattice_from_category}, $L_x$ is a Garside lattice, and one remarks that the edge relation on $\Cay(G_x,S)$ coincides with the one given in Section~\ref{sec:lattice_no_quotient}, and the edge relation on the quotient $\Cay(G_x,S) / \<\Delta_x\>$ coincides with the one given in Section~\ref{sec:quotient}.
\ep

\subsection{Artin complexes}

\bthm \label{thm:Artin complex}
Let $X$ be the Artin complex of the Artin-Tits group of type $\widetilde A_{n-1}$ (cf. Section~\ref{subsec:Artin}). Then the 1-skeleton of $X$ is a weakly modular graph.
\ethm

\begin{proof}
Let $A_i$ be as defined in Section~\ref{subsec:Artin}. We represent vertices of $X$ by left cosets of the form $gA_i$.
Let $P=X^{(0)}\times \Z$. 
We put a weak order on $P$ as follows. Define $(g_1A_{i_1},n)\le (g_2A_{i_2},m)$ if one of the following holds:
\begin{enumerate}
	\item $n=m$, $g_1A_{i_1}\cap g_2A_{i_2}\neq\emptyset$ and $i_1\le i_2$;
	\item $m=n+1$, $g_1A_{i_1}\cap g_2A_{i_2}\neq\emptyset$ and $i_2\le i_1$.
\end{enumerate}
Note that $(P,\le)$ is a weakly ordered set. If $(x,n)\le (y,m)$, then $x$ and $y$ are adjacent in $X$. Conversely, if $x$ and $y$ are adjacent, then either $(x,n)\ge (y,n)$ or $(x,n)\ge (y,n-1)$. As $X^{(1)}$ is connected, we deduce that for any $(x,n)\in P$ and any $y\in X$, there exists $m\in\Z$ such that we can find a weak chain from $(y,m)$ to $(x,n)$.
The weakly ordered set $P$ is homogeneous, indeed, for $a\le b$ in $P$, we can define the length function $\ell(a\le b)$ to be the maximal length of chains in $P$ from $a$ to $b$. 

Let $\varphi: P\to P$ be the map sending $(x,n)$ to $(x,n+1)$. Then $\varphi$ is an automorphism of weakly ordered set. Then $\varphi$ generates the weak order $\le$. Theorem~\ref{thm:lattice ltg} (2) holds true as $X_{\varphi}$ is homeomorphic to $X\times \mathbb R$ and $X$ is simply-connected. Theorem~\ref{thm:lattice ltg} (3) holds true by \cite[Theorem 4.3]{haettel2021lattices}. Thus Theorem~\ref{thm:lattice ltg} implies $\le$ generates a partial order $\le_t$ on $P$. The previous paragraph implies any two elements in $(P,\le_t)$ have lower bound, thus $(P,\le_t)$ is a homogeneous lattice by Theorem~\ref{thm:lattice ltg}. As $(P,\le_t)$ is homogeneous, we know any lower bounded subset in $(P,\le_t)$ has a meet and any upper bounded subset in $(P,\le_t)$ has a join. Now the automorphism $\varphi$ gives an action of $\Z$ on $(P,\le_t)$. The previous paragraph also implies for any $a,b\in (P,\le_t)$, there exists $k\in \Z$ such that $-k\cdot a\le b\le k\cdot a$.
Note that $X^{(1)}$ coincides with the quotient graph of $(P,\le_t)$ by $\Z$ as defined in Section~\ref{sec:quotient}. Thus we are done by Theorem~\ref{thm:quotient wm}.
\end{proof}

\subsection{Arc complexes}

The following is nothing more than a topological description of the Artin complex of the Artin-Tits group of type $\widetilde A_{n-1}$. Let $\Sigma$ denote a $2$-sphere with $n+2$ punctures $\{N,S,p_1,\dots,p_n\}$ with two distinguished punctures $N,S$ which could thought of as the North pole and the South pole of $\Sigma$. The punctures $p_1,\dots,p_n$ may be thought as cyclically ordered on the equator of the $2$-sphere.

\begin{figure}[H]
\centering
\includegraphics[width=5cm]{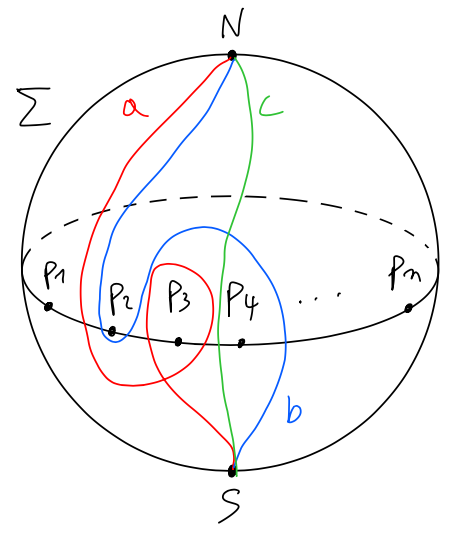}
\caption{Arcs on the punctured sphere $\Sigma$: $a$ is adjacent to $b$ and $c$}
\label{fig:arc_complex}
\end{figure}

Let ${\cal A}(\Sigma)$ denote the following simplicial complex. Its vertex set consists of isotopy classes of arcs in $\Sigma$ from $N$ to $S$. Two vertices are adjacent if they can be realized disjointly. Then ${\cal A}(\Sigma)$ is the associated flag simplicial complex, see Figure~\ref{fig:arc_complex}. According to~\cite[Lemma~2.5]{wahl}, this arc complex ${\cal A}(\Sigma)$ is contractible.

\bthm \label{thm:arc_complex_wm}
The arc complex ${\cal A}(\Sigma)$ is a weakly modular graph.
\ethm

\begin{figure}[H]
\centering\includegraphics[width=12cm]{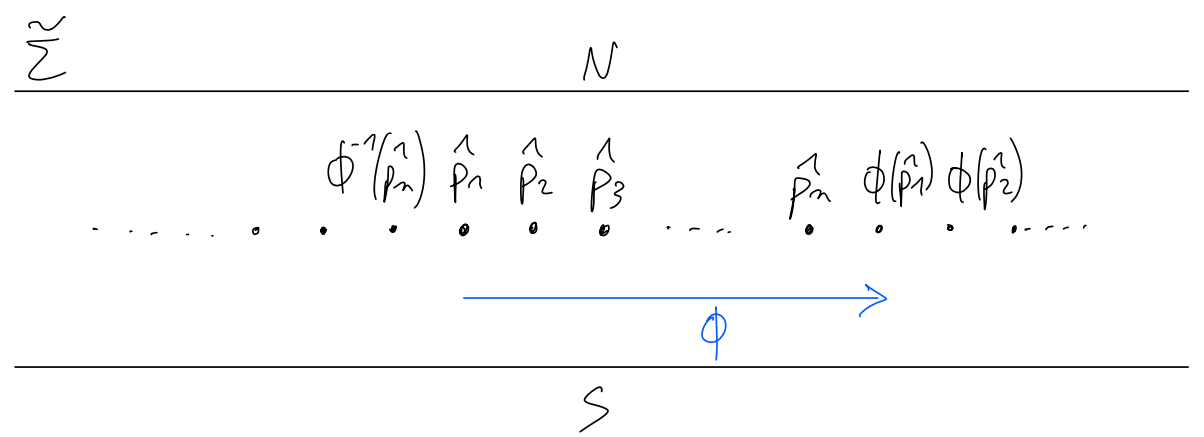}
\caption{The cyclic cover $\tilde{\Sigma}$.}
\label{fig:cyclic_cover}
\end{figure}

\bp
Let us consider the cyclic cover $\tilde{\Sigma}$ of $\Sigma$ over the set of poles $\{N,S\}$. More precisely, consider the group morphism $\psi:\pi_1(\Sigma) \ra \Z$ sending a loop around $N$ to $+1$, a loop around $S$ to $-1$, and a loop around any $p_i$ to $0$. The cyclic covering $\tilde{\Sigma}$ of $\Sigma$ associated to $\Ker \psi$ is homeomorphic to a $2$-sphere, with $2$ distinguished punctures denoted $N,S$, and infinitely many punctures forming $n$ $\phi$-orbits $\{\phi^k(\hat{p_i}), 1 \leq i \leq n, k \in \Z\}$, where $\phi$ is a generator the of deck transformation group of $\tilde{\Sigma}$ and $\hat{p_i}$ is a lift of $p_i$ in $\tilde{\Sigma}$, see Figure~\ref{fig:cyclic_cover}.

Let us consider the set $P$ of isotopy classes of arcs $a$ from $N$ to $S$ in $\tilde{\Sigma}$, which are lifts of arcs in ${\cal A}(\Sigma)$. There is an induced induced action of $\phi$ on $P$.

We will put a weak order on $P$ as follows. Say that $a \leq a'$ if $a'=a$, $a'=\phi(a)$, or $a'$ is disjoint from $\phi(a)$, and $a'$ separates $a$ and $\phi(a)$. Since the arc complex ${\cal A}(\Sigma)$ is connected, we deduce that $\phi : P \ra P$ generates $\leq$.

\begin{figure}[H]
\centering
\includegraphics[width=5cm]{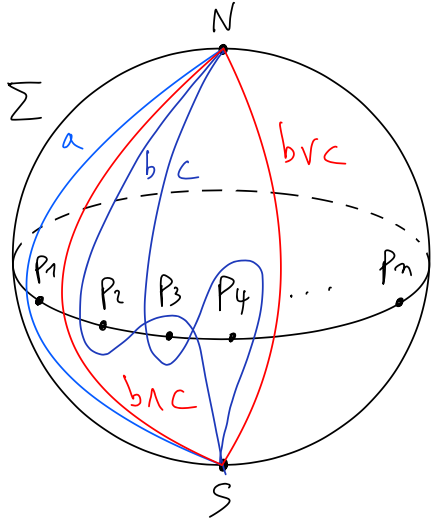}
\caption{The lattice property: the meet $b \wedge c$ and the join $b \vee c$ of the two arcs $b,c$ in the interval $[a,\phi(a)]$.}
\label{fig:lattice}
\end{figure}

We will show that we can apply Theorem~\ref{thm:lattice ltg}.
\ben
\item By definition of the weak order, for any $a \in P$, we have $\phi(a) > a$.
\item Since the arc complex ${\cal A}(\Sigma)$ is simply connected, we deduce that $X_\phi$ is simply connected.
\item For any $a \in P$, the interval $[a,\phi(a)]$ is isomorphic to the lattice of cut-curves, see~\cite{bessis} and \cite{haettel2021lattices} for the proof of the lattice property due to Crisp and McCammond (unpublished). The lattice property is quite geometric in this case: roughly speaking, the meet of two arcs $b,c$ in the interval $[a,\phi(a)]$ may be defined as the "westmost" part of $b \cup c$, see Figure~\ref{fig:lattice}.
\een

According to Theorem~\ref{thm:lattice ltg}, we deduce that $(P,\leq_t)$ is a lattice. Moreover, given any $a,a' \in P$, there exists $k \in \N$ such that $a \leq_t \phi^k(a')$. So we can apply Theorem~\ref{thm:quotient wm}, and deduce that the quotient graph, which coincides with ${\cal A}(\Sigma)$, is a weakly modular graph.
\ep

It turns out that this arc complex coincides with the above Artin complex.

\bpro 
\label{prop:iso}
The arc complex ${\cal A}(\Sigma)$ is isomorphic to the Artin complex of the affine Artin-Tits group of type $A(\widetilde A_{n-1})$. \epro

\bp Note that the mapping class group $G=\Mod(\Sigma, \{N,S\})$ fixing the set of North and South poles act by simplicial automorphisms on ${\cal A}(\Sigma)$. According to~\cite{charney_crisp_Artin_MCG}, $G$ is isomorphic to the semidirect product $A(\widetilde A_{n-1}) \rtimes \Z/n\Z$, where $A(\widetilde A_{n-1})$ is the affine Artin group of type $\widetilde A_{n-1}$ and $\Z/n\Z$ acts by rotations on the defining graph of $A(\widetilde A_{n-1})$.

\mk

\begin{figure}[H]
\centering
\includegraphics[width=5cm]{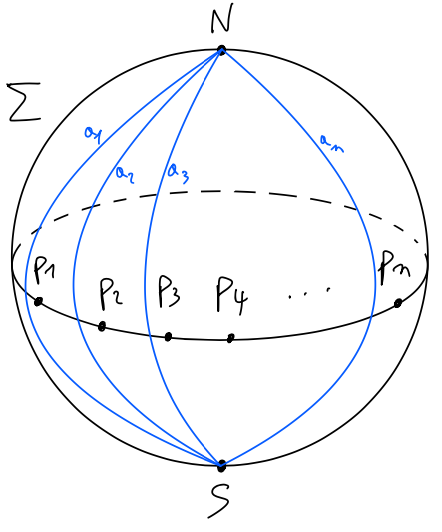}
\caption{The fundamental simplex $a_1,a_2,\dots,a_n$ of the arc complex.}
\label{fig:simplex}
\end{figure}

Note that the subgroup $A=A(\widetilde A_{n-1})$ of $G$ acts on ${\cal A}(\Sigma)$, with strict fundamental domain the $n$-simplex consisting of the $n$ meridians $a_1,\dots,a_n$ "separating" the points $p_1,\dots,p_n$ in this cyclic order. More precisely, $a_i$ crosses the equatorial arc between $p_i$ and $p_{i+1}$, see Figure~\ref{fig:simplex}.

Note that the stabilizer of each $a_i$ is naturally isomorphic to the $n$-strand braid group, generated by $n-1$ consecutive standard generators of $A(\widetilde A_{n-1})$. Similarly, the stabilizer of $(a_{i_1},\dots,a_{i_k})$ is isomorphic to the corresponding product of braid groups. We deduce that both ${\cal A}(\Sigma)$ and the Artin complex of affine type $A(\widetilde A_{n-1})$ are universal covers of the simplex of groups defining the Artin complex of affine type $A(\widetilde A_{n-1})$.
\ep

\addcontentsline{toc}{section}{References}

\bibliographystyle{alpha}
\bibliography{bibli}

\newcommand{\etalchar}[1]{$^{#1}$}
\def\polhk#1{\setbox0=\hbox{#1}{\ooalign{\hidewidth
  \lower1.5ex\hbox{`}\hidewidth\crcr\unhbox0}}}
\begin{thebibliography}{CMSZ93b}

\bibitem[AB08]{abramenko_brown}
Peter Abramenko and Kenneth~S. Brown.
\newblock {\em {Buildings. Theory and Applications.\!}}
\newblock Grad. Text. Math. Springer, 2008.

\bibitem[BC96]{bandelt1996helly}
Hans-J{\"u}rgen Bandelt and Victor Chepoi.
\newblock A helly theorem in weakly modular space.
\newblock {\em Discrete Mathematics}, 160(1-3):25--39, 1996.

\bibitem[BC06]{bessis2006non}
David Bessis and Ruth Corran.
\newblock Non-crossing partitions of type (e, e, r).
\newblock {\em Advances in Mathematics}, 202(1):1--49, 2006.

\bibitem[BCC{\etalchar{+}}13]{brevsar2013bucolic}
Bo{\v{s}}tjan Bre{\v{s}}ar, J{\'e}r{\'e}mie Chalopin, Victor Chepoi, Tanja
  Gologranc, and Damian Osajda.
\newblock Bucolic complexes.
\newblock {\em Advances in mathematics}, 243:127--167, 2013.

\bibitem[BCL16]{bader_caprace_lecureux}
Uri Bader, Pierre-Emmanuel Caprace, and Jean L{\~A}{\copyright}cureux.
\newblock {On the linearity of lattices in affine buildings and ergodicity of
  the singular Cartan flow}.
\newblock 2016.
\newblock {arXiv:1608.06265}.

\bibitem[Bes99]{bestvina_artin}
Mladen Bestvina.
\newblock Non-positively curved aspects of {A}rtin groups of finite type.
\newblock {\em Geom. Topol.}, 3:269--302, 1999.

\bibitem[Bes03]{bessis}
David Bessis.
\newblock The dual braid monoid.
\newblock {\em Ann. Sci. \'Ecole Norm. Sup. (4)}, 36(5):647--683, 2003.

\bibitem[Bes06]{bessis2006garside}
David Bessis.
\newblock Garside categories, periodic loops and cyclic sets.
\newblock {\em arXiv preprint math/0610778}, 2006.

\bibitem[Bes15]{bessis2015finite}
David Bessis.
\newblock {Finite complex reflection arrangements are K ($\pi$, 1)}.
\newblock {\em Annals of mathematics}, pages 809--904, 2015.

\bibitem[BF14]{bader_furman_icm}
Uri Bader and Alex Furman.
\newblock Boundaries, rigidity of representations, and {L}yapunov exponents.
\newblock In {\em Proceedings of the {I}nternational {C}ongress of
  {M}athematicians---{S}eoul 2014. {V}ol. {III}}, pages 71--96. Kyung Moon Sa,
  Seoul, 2014.

\bibitem[Blu21]{blufstein2021parabolic}
Martin~Axel Blufstein.
\newblock {Parabolic subgroups of two-dimensional Artin groups and
  systolic-by-function complexes}.
\newblock {\em arXiv preprint arXiv:2108.04929}, 2021.

\bibitem[CC05]{charney_crisp_Artin_MCG}
Ruth Charney and John Crisp.
\newblock Automorphism groups of some affine and finite type {A}rtin groups.
\newblock {\em Math. Res. Lett.}, 12(2-3):321--333, 2005.

\bibitem[CCG{\etalchar{+}}20]{chalopin2020helly}
J{\'e}r{\'e}mie Chalopin, Victor Chepoi, Anthony Genevois, Hiroshi Hirai, and
  Damian Osajda.
\newblock Helly groups.
\newblock {\em arXiv preprint arXiv:2002.06895}, 2020.

\bibitem[CCHO21]{chalopin_chepoi_hirai_osajda}
J{\'e}r{\'e}mie Chalopin, Victor Chepoi, Hiroshi Hirai, and Damian Osajda.
\newblock Weakly modular graphs and nonpositive curvature.
\newblock {\em Mem. Amer. Math. Soc.}, 2021.

\bibitem[CD95]{charney_davis_kpi1}
Ruth Charney and Michael~W. Davis.
\newblock The {$K(\pi,1)$}-problem for hyperplane complements associated to
  infinite reflection groups.
\newblock {\em J. Amer. Math. Soc.}, 8(3):597--627, 1995.

\bibitem[Che89]{chepoi1989classification}
VD~Chepo{\i}.
\newblock Classification of graphs by means of metric triangles.
\newblock {\em Metody Diskret. Analiz}, 49:75--93, 1989.

\bibitem[Che00]{chepoi2000graphs}
Victor Chepoi.
\newblock Graphs of some cat (0) complexes.
\newblock {\em Advances in Applied Mathematics}, 24(2):125--179, 2000.

\bibitem[CLL15]{corran2015braid}
Ruth Corran, Eon-Kyung Lee, and Sang-Jin Lee.
\newblock Braid groups of imprimitive complex reflection groups.
\newblock {\em Journal of Algebra}, 427:387--425, 2015.

\bibitem[CMSZ93a]{cartwright_mantero_steger_zappa_I}
Donald~I. Cartwright, Anna~Maria Mantero, Tim Steger, and Anna Zappa.
\newblock Groups acting simply transitively on the vertices of a building of
  type {$\tilde A_2$}. {I}.
\newblock {\em Geom. Dedicata}, 47(2):143--166, 1993.

\bibitem[CMSZ93b]{cartwright_mantero_steger_zappa_II}
Donald~I. Cartwright, Anna~Maria Mantero, Tim Steger, and Anna Zappa.
\newblock Groups acting simply transitively on the vertices of a building of
  type {$\tilde A_2$}. {II}. {T}he cases {$q=2$} and {$q=3$}.
\newblock {\em Geom. Dedicata}, 47(2):167--223, 1993.

\bibitem[CMV20]{cumplido2020parabolic}
Mar{\'\i}a Cumplido, Alexandre Martin, and Nicolas Vaskou.
\newblock {Parabolic subgroups of large-type Artin groups}.
\newblock {\em arXiv preprint arXiv:2012.02693}, 2020.

\bibitem[CP05]{crisp2005representations}
John Crisp and Luis Paris.
\newblock {Representations of the braid group by automorphisms of groups,
  invariants of links, and Garside groups}.
\newblock {\em Pacific journal of mathematics}, 221(1):1--27, 2005.

\bibitem[CP11]{corran2011new}
Ruth Corran and Matthieu Picantin.
\newblock A new garside structure for the braid groups of type (e, e, r).
\newblock {\em Journal of the London Mathematical Society}, 84(3):689--711,
  2011.

\bibitem[CW17a]{calvez_wiest_2}
Matthieu Calvez and Bert Wiest.
\newblock Acylindrical hyperbolicity and {A}rtin-{T}its groups of spherical
  type.
\newblock {\em Geom. Dedicata}, 191:199--215, 2017.

\bibitem[CW17b]{calvez_wiest_1}
Matthieu Calvez and Bert Wiest.
\newblock Curve graphs and {G}arside groups.
\newblock {\em Geom. Dedicata}, 188:195--213, 2017.

\bibitem[Deh15]{garside}
Patrick Dehornoy.
\newblock {\em Foundations of {G}arside theory}, volume~22 of {\em EMS Tracts
  in Mathematics}.
\newblock European Mathematical Society (EMS), Z\"urich, 2015.
\newblock With Fran\c{c}ois Digne, Eddy Godelle, Daan Krammer and Jean Michel,
  Contributor name on title page: Daan Kramer.

\bibitem[Del72]{deligne1972immeubles}
Pierre Deligne.
\newblock Les immeubles des groupes de tresses g{\'e}n{\'e}ralis{\'e}s.
\newblock {\em Inventiones mathematicae}, 17(4):273--302, 1972.

\bibitem[Dig06]{digne2006presentations}
Fran{\c{c}}ois Digne.
\newblock Pr{\'e}sentations duales des groupes de tresses de type affine a.
\newblock {\em Commentarii Mathematici Helvetici}, 81(1):23--47, 2006.

\bibitem[Dig12]{digne2012garside}
Fran{\c{c}}ois Digne.
\newblock A garside presentation for artin-tits groups of type $\widetilde
  c_n$.
\newblock In {\em Annales de l'Institut Fourier}, volume~62, pages 641--666,
  2012.

\bibitem[dlSLW23]{lecureux_delasalle_witzel_T}
Mikael de~la Salle, Jean L\'ecureux, and Stefan Witzel.
\newblock {Strong property (T), weak amenability and $\ell^p$-cohomology in
  $\tilde{A}_2$-buildings}.
\newblock {\em Ann. Sci. \'ecole Norm. Sup.}, 2023.

\bibitem[Gen17]{genevois2017cubical}
Anthony Genevois.
\newblock Cubical-like geometry of quasi-median graphs and applications to
  geometric group theory.
\newblock {\em arXiv preprint arXiv:1712.01618}, 2017.

\bibitem[Gro87]{gromov_hyperbolic_groups}
M.~Gromov.
\newblock Hyperbolic groups.
\newblock In {\em Essays in group theory}, volume~8 of {\em Math. Sci. Res.
  Inst. Publ.}, pages 75--263. Springer, New York, 1987.

\bibitem[Hae20]{haettel_hyprig}
Thomas Haettel.
\newblock Hyperbolic rigidity of higher rank lattices.
\newblock {\em with an appendix by Vincent Guirardel et Camille Horbez, Ann.
  Sci. \'{E}c. Norm. Sup\'{e}r. (2)}, 53:439--468, 2020.

\bibitem[Hae21a]{haettel_injective_buildings}
Thomas Haettel.
\newblock {Injective metrics on buildings and symmetric spaces}.
\newblock {arXiv:2101.09367}, 2021.

\bibitem[Hae21b]{haettel2021lattices}
Thomas Haettel.
\newblock Lattices, injective metrics and the $k(\pi,1)$ conjecture.
\newblock {\em arXiv preprint arXiv:2109.07891}, 2021.

\bibitem[Hir20]{hirai_uniform_modular}
Hiroshi Hirai.
\newblock Uniform modular lattices and affine buildings.
\newblock {\em Adv. Geom.}, 20(3):375--390, 2020.

\bibitem[HO19]{huang_osajda_helly}
Jingyin Huang and Damian Osajda.
\newblock {Helly meets Garside and Artin}.
\newblock {\em arXiv:1904.09060}, 2019.

\bibitem[HO20]{huang2020large}
Jingyin Huang and Damian Osajda.
\newblock Large-type artin groups are systolic.
\newblock {\em Proceedings of the London Mathematical Society}, 120(1):95--123,
  2020.

\bibitem[HO21]{haettel2021locally}
Thomas Haettel and Damian Osajda.
\newblock Locally elliptic actions, torsion groups, and nonpositively curved
  spaces.
\newblock {\em arXiv preprint arXiv:2110.12431}, 2021.

\bibitem[Hod20]{hoda:crystallographic}
Nima Hoda.
\newblock Crystallographic {H}elly groups.
\newblock {\em arXiv preprint arXiv:2010.07407}, 2020.

\bibitem[JS06]{januszkiewicz2006simplicial}
Tadeusz Januszkiewicz and Jacek Swiatkowski.
\newblock Simplicial nonpositive curvature.
\newblock {\em Publications Math{\'e}matiques de l'Institut des Hautes
  {\'E}tudes Scientifiques}, 104(1):1--85, 2006.

\bibitem[Lan13]{lang}
Urs Lang.
\newblock Injective hulls of certain discrete metric spaces and groups.
\newblock {\em J. Topol. Anal.}, 5(3):297--331, 2013.

\bibitem[McC15]{mccammond2015dual}
Jon McCammond.
\newblock Dual euclidean artin groups and the failure of the lattice property.
\newblock {\em Journal of Algebra}, 437:308--343, 2015.

\bibitem[MS17]{mccammond2017artin}
Jon McCammond and Robert Sulway.
\newblock Artin groups of euclidean type.
\newblock {\em Inventiones mathematicae}, 210(1):231--282, 2017.

\bibitem[Mun19]{munro2019weak}
Zachary Munro.
\newblock Weak modularity and $\widetilde {A}_n$ buildings.
\newblock {\em arXiv preprint arXiv:1906.10259}, 2019.

\bibitem[PS16]{przytycki_schwer_systolic_buildings}
Piotr Przytycki and Petra Schwer.
\newblock Systolizing buildings.
\newblock {\em Groups Geom. Dyn.}, 10(1):241--277, 2016.

\bibitem[Rad17]{radu_exotic_lattice}
Nicolas Radu.
\newblock A lattice in a residually non-{D}esarguesian {$\tilde A_2$}-building.
\newblock {\em Bull. Lond. Math. Soc.}, 49(2):274--290, 2017.

\bibitem[Soe21]{soergel2021systolic}
Mireille Soergel.
\newblock Systolic complexes and group presentations.
\newblock {\em arXiv preprint arXiv:2105.01345}, 2021.

\bibitem[Tit66]{tits_artin}
J.~Tits.
\newblock Normalisateurs de tores. {I}. {G}roupes de {C}oxeter \'{e}tendus.
\newblock {\em J. Algebra}, 4:96--116, 1966.

\bibitem[VdL83]{vanderlek}
Harm Van~der Lek.
\newblock {The homotopy type of complex hyperplane complements}.
\newblock {1983}.
\newblock {PhD thesis, Radboud University of Nijmegen}.

\bibitem[Wah08]{wahl}
Nathalie Wahl.
\newblock Homological stability for the mapping class groups of non-orientable
  surfaces.
\newblock {\em Invent. Math.}, 171(2):389--424, 2008.

\end{thebibliography}

\end{document}